\documentclass[aos,preprint]{imsart}

\RequirePackage{amsthm,amsmath}
\RequirePackage[numbers]{natbib}%
\RequirePackage[colorlinks,citecolor=blue,urlcolor=blue]{hyperref}

\usepackage{amsgen,amsmath,amstext,amsbsy,amsopn,amssymb,latexsym,comment}
\usepackage{array}
\usepackage{cite}
\usepackage{multirow}
\usepackage{graphicx}
\usepackage{amssymb}
\usepackage{bibunits}
\usepackage[dvipsnames,table,dvipsnames*, svgnames*, hyperref]{xcolor}
\arxiv{arXiv:0000.0000}

\newcommand{\CIZR}{{\rm CI}_{\alpha}\left(\widehat{\beta}, \ell_{q},Z\right)}
\newcommand{\gCIZR}{{\rm CI}_{\alpha}\left(\widehat{\beta},\ell_{q}, Z\right)}

\newcommand{\lossnull}{\|\widehat{\beta}(z)-\betanull\|_{q}^2}
\newcommand{\lossnullr}{\|\widehat{\beta}(Z)-\betanull\|_{q}^2}
\newcommand{\lest}{\widehat{L}_q}
\newcommand{\gloss}{\|\widehat{\beta}-\beta(\theta)\|_q^2}
\newcommand{\glosss}{\|\widehat{\beta}-\beta\|_q^2}

\newcommand{\liminfnp}{\liminf\limits_{n,p \rightarrow \infty}}
\newcommand{\limsupnp}{\limsup\limits_{n,p \rightarrow \infty}}

\newtheorem{Proposition}{Proposition}
\newtheorem{Lemma}{Lemma}

\newtheorem{Theorem}{Theorem}

\newtheorem{Remark}{Remark}

\newcommand{\R}{\mathbb{R}}
\newcommand{\E}{{\mathbf{E}}}
\newcommand{\FF}{{\mathcal{F}}}

\newcommand{\II}{{\mathcal{I}}}

\newcommand{\PP}{{\mathbb{P}}}

\newcommand{\RR}{\mathbf{L}}

\newcommand{\Xj}{X_{\cdot j}}
\newcommand{\argmin}{\arg\min}
\newcommand{\Mv}{M_2}
\newcommand{\Md}{M_1}
\newcommand{\TV}{L_{1}}

\newcommand{\pjoint}{\PP_{Z,\theta \sim \pi}}
\newcommand{\betanull}{\beta^{*}}

\newcommand{\Thetas}{\Theta_1}
\newcommand{\Thetal}{\Theta_2}

\newcommand{\suppx}{{\rm supp}(x)}

\newcommand{\regimea}{k_1\leq k_2\lesssim \frac{\sqrt{n}}{\log p}}
\newcommand{\regimeb}{k_1\lesssim \frac{\sqrt{n}}{\log p} \lesssim k_2 \lesssim \frac{n}{\log p}}
\newcommand{\regimec}{\frac{\sqrt{n}}{\log p}\lesssim k_1\leq  k_2 \lesssim \frac{n}{\log p}}

\newcommand{\nf}{n_1}
\newcommand{\ns}{n_2}
\newcommand{\xa}{X^{\left(1\right)}}

\newcommand{\ya}{y^{\left(1\right)}}

\newcommand{\distalter}{d}

\newcommand{\disttwonorm}{{\rm dist}}
\newcommand{\betal}{\widehat{\beta}^{L}}
\newcommand{\betasl}{\widehat{\beta}^{SL}}

\startlocaldefs
\numberwithin{equation}{section}
\theoremstyle{plain}

\endlocaldefs

\defaultbibliographystyle{plain}
\defaultbibliography{HDRef}

\begin{document}
\begin{frontmatter}
\title{Accuracy Assessment for High-dimensional Linear Regression\thanksref{T1}}
\runtitle{High-dimensional Accuracy Assessment}
\thankstext{T1}{The research was supported in part by NSF Grants DMS-1208982 and DMS-1403708, and NIH Grant R01 CA127334.}

\begin{aug}
\author{\fnms{T. Tony} \snm{Cai}\ead[label=e1]{tcai@wharton.upenn.edu}},
\and
\author{\fnms{Zijian} \snm{Guo}\ead[label=e2]{zijguo@wharton.upenn.edu}
\ead[label=u1,url]{URL: http://www-stat.wharton.upenn.edu/$\sim$tcai/}}
\runauthor{T. T. Cai and Z. Guo}
\affiliation{University of Pennsylvania}
\address{DEPARTMENT OF STATISTICS\\
THE WHARTON SCHOOL\\
UNIVERSITY OF PENNSYLVANIA\\
PHILADELPHIA, PENNSYLVANIA 19104\\
USA\\
\printead{e1}\\
\phantom{E-mail:\ }\printead*{e2}\\
\printead*{u1}\phantom{URL:\ }\\
}

\end{aug}


\begin{abstract}
This paper considers point and interval estimation of  the $\ell_q$ loss of an estimator in high-dimensional linear regression with random design. We establish the minimax rate for estimating the $\ell_{q}$ loss and the minimax expected length of confidence intervals for the $\ell_{q}$ loss of rate-optimal estimators of the regression vector, including commonly used estimators such as Lasso, scaled Lasso, square-root Lasso and Dantzig Selector. Adaptivity  of the confidence intervals for the $\ell_{q}$ loss is also studied. Both the setting of known identity design covariance matrix and known noise level and the setting of unknown design covariance matrix and unknown noise level are studied.  The results reveal interesting and significant differences between estimating the $\ell_2$ loss and $\ell_q$ loss with $1\le q <2$ as well as between the two settings.

New technical tools are developed to establish rate sharp lower bounds for the minimax estimation error and  the expected length of minimax and adaptive confidence intervals for the $\ell_q$ loss. A significant difference between loss estimation and the traditional parameter estimation is that for loss estimation  the constraint is on the performance of the estimator of the regression vector, but the lower bounds are on the difficulty of estimating its $\ell_q$ loss. 
The technical tools developed in this paper can also be of independent interest.
\end{abstract}

\begin{keyword}[class=MSC]
\kwd[Primary ]{62G15}
\kwd[; secondary ]{62C20}\kwd{62H35}
\end{keyword}

\begin{keyword}
\kwd{Accuracy assessment, adaptivity, confidence interval, high-dimensional linear regression, loss estimation, minimax lower bound, minimaxity,  sparsity.}
\end{keyword}
\end{frontmatter}

\begin{bibunit}

\section{Introduction}

In many applications, the goal of statistical inference is not only to construct a good estimator, but also to provide a measure of accuracy for this estimator. In classical statistics, when the parameter of interest is one-dimensional, this is achieved in the form of a standard error or a confidence interval. A prototypical example is the inference for a binomial proportion, where often not only an estimate of the proportion but also its margin of error are given. Accuracy measures of an estimation procedure have also been used as a tool for the empirical selection of tuning parameters. A well known example is Stein's Unbiased Risk Estimate (SURE), which has been an effective tool for the construction of data-driven adaptive estimators in normal means estimation,  nonparametric signal recovery, covariance matrix estimation, and other problems. See, for instance, \citep{stein1981estimation,li1985stein,donoho1995adapting,cai2009data,yi2013sure}. The commonly used cross-validation methods  can also be viewed as a useful tool based on the idea of empirical assessment of accuracy.

In this paper, we consider the problem of estimating the loss of a given estimator in the setting of high-dimensional linear regression, where one observes $(X, y)$ with $X \in \R^{n\times p}$ and $y \in \R^{n}$, and for $1\leq i\leq n,$
\[
y_i=X_{i\cdot}\beta+\epsilon_i.
\]
Here $\beta \in \R^{p}$ is the regression vector,  $X_{i\cdot} \stackrel{iid}{\sim} N_p(0,\Sigma)$ are the rows of $X$, and the errors $\epsilon_i \stackrel{iid}{\sim} N(0,\sigma^2)$ are independent of $X$.  This high-dimensional linear model has been well studied in the literature, with the main focus on estimation of $\beta$. Several penalized/constrained $\ell_1$ minimization methods, including Lasso \citep{tibshirani1996regression}, Dantzig selector \citep{candes2007dantzig}, scaled Lasso \citep{sun2012scaled} and square-root Lasso \citep{belloni2011square}, have been proposed.  These methods have been shown to work well in applications and produce interpretable estimates of $\beta$ when $\beta$ is assumed to be sparse. Theoretically, with a properly chosen tuning parameter, these estimators  achieve the optimal rate of convergence over collections of sparse parameter spaces. See, for example, \citep{candes2007dantzig,sun2012scaled,belloni2011square,raskutti2011minimax,bickel2009simultaneous,buhlmann2011statistics,verzelen2012minimax}.

For a given estimator $\widehat{\beta}$, the $\ell_q$ loss $\|\widehat{\beta}-\beta\|_q^2$ with $1\le q\le 2$ is commonly used as a metric of accuracy for $\widehat{\beta}$.  We consider in the present paper both point and interval estimation of  the $\ell_q$ loss $\|\widehat{\beta}-\beta\|_q^2$ for a given $\widehat{\beta}$.
Note that the loss $\|\widehat{\beta}-\beta\|_q^2$ is a random quantity, depending on  both the estimator $\widehat{\beta}$ and the parameter $\beta$. For such a random quantity, prediction and prediction interval are ususally used for point and interval estimation, respectively. However, we slightly abuse the terminologies in the present paper by using estimation and confidence interval to represent the point and interval estimators of the loss $\|\widehat{\beta}-\beta\|_q^2$. Since the $\ell_{q}$ loss depends on the estimator $\widehat{\beta}$, it is necessary to specify the estimator in the discussion of loss estimation. Throughout this paper, we restrict our attention to a broad collection of estimators $\widehat{\beta}$ that perform well at least at one interior point or a small subset of the parameter space. This collection of estimators includes most state-of-art estimators such as Lasso, Dantzig selector, scaled Lasso and square-root Lasso.

High-dimensional linear regression has been well studied in two settings. One is the setting with known design covariance matrix $\Sigma={\rm I}$ and known noise level $\sigma=\sigma_0$ and sparse $\beta$. See for example, \citep{donoho2011noise,bayati2012lasso,nickl2013confidence,verzelen2012minimax,thrampoulidis2015asymptotically,janson2015eigenprism,cai2015regci,arias2011global,ingster2010detection}. Another commonly considered setting is sparse $\beta$ with unknown $\Sigma$ and $\sigma$. We study point and interval estimation of the $\ell_q$ loss $\|\widehat{\beta}-\beta\|_q^2$ in both settings. Specifically, we consider the parameter space $\Theta_0(k)$ introduced in \eqref{eq: parameter space with fixed variance and known design}, which consists of $k$-sparse signals $\beta$ with known design covariance matrix $\Sigma={\rm I}$ and known noise level $\sigma=\sigma_0$, and $\Theta(k)$ defined in \eqref{eq: parameter space with unknown variance}, which consists of $k$-sparse signals with unknown $\Sigma$ and $\sigma$.

\subsection{Our contributions}


The present paper studies the minimax and adaptive estimation of the loss $\|\widehat{\beta}-\beta\|_q^2$ for a given estimator $\widehat{\beta}$ and  the minimax expected length and adaptivity of confidence intervals for the loss.
 A major step in our analysis is to establish rate sharp lower bounds for the minimax estimation error  and the minimax expected length of confidence intervals for the $\ell_{q}$ loss over $\Theta_0(k)$ and $\Theta(k)$ for a broad class of estimators of $\beta$, which contains the subclass of rate-optimal estimators. We then focus on the estimation of the loss of rate-optimal estimators and take the Lasso and scaled Lasso estimators as generic examples. For these rate-optimal estimators, we propose procedures for point estimation as well as confidence intervals for their $\ell_q$ losses. It is shown that the proposed procedures achieve the corresponding lower bounds up to a constant factor.   These results together establish the minimax rates for estimating the $\ell_{q}$ loss of rate-optimal estimators over $\Theta_0(k)$ and $\Theta(k)$.  The analysis shows interesting and significant differences between estimating the $\ell_2$ loss and $\ell_q$ loss with $1\le q <2$ as well as between the two parameter spaces $\Theta(k)$ and $\Theta_0(k)$.
\begin{itemize}
\item The minimax rate for estimating $\|\widehat{\beta}-\beta\|_2^2$ over $\Theta_0(k)$ is $\min\left\{\frac{1}{\sqrt{n}}, k\frac{\log p}{n}\right\}$ and over $\Theta(k)$ is $k \frac{\log p}{n}$. So loss estimation is much easier with the prior information $\Sigma={\rm I}$ and $\sigma=\sigma_0$ when ${\sqrt{n} \over \log p}\ll k \lesssim {n\over \log p}$.

\item The minimax rate for estimating $\|\widehat{\beta}-\beta\|_q^2$ with $1\leq q< 2$ over both $\Theta_0(k)$ and $\Theta(k)$ is $k^{\frac{2}{q}}\frac{\log p}{n}$.
\end{itemize}

In the regime $\frac{\sqrt{n}}{\log p}\ll k \lesssim \frac{n}{\log p},$ a practical loss estimator is proposed for estimating the $\ell_2$ loss and shown to achieve the optimal convergence rate $\frac{1}{\sqrt{n}}$  adaptively over $\Theta_0(k)$. We say {\it estimation of loss is  impossible} if the minimax rate can be achieved by the trivial estimator 0, which means that the estimation accuracy of the loss is at least of the same order as the loss itself.  In all other considered cases, estimation of loss is shown to be impossible. These results indicate that loss estimation is  difficult.

We then turn to the construction of confidence intervals for the $\ell_q$ loss. A confidence interval for the loss is useful even when it is ``impossible" to estimate the loss, as a confidence interval can provide non-trivial upper and lower bounds for the loss. In terms of convergence rate over $\Theta_0(k)$ or $\Theta(k)$, the minimax rate of the expected length of confidence intervals for the $\ell_q$ loss, $\glosss$, of any rate-optimal estimator $\widehat{\beta}$ coincides with the minimax estimation rate. We also consider the adaptivity of confidence intervals for the $\ell_q$ loss of any rate-optimal estimator $\widehat{\beta}$. (The framework for adaptive confidence intervals is discussed in detail in Section \ref{framework.sec}.)
Regarding confidence intervals for the $\ell_2$ loss in the case of known $\Sigma={\rm I}$ and $\sigma=\sigma_0$,  a procedure is proposed  and is shown to achieve the optimal length $\frac{1}{\sqrt{n}}$ adaptively over $\Theta_0(k)$ for $\frac{\sqrt{n}}{\log p}\lesssim k  \lesssim \frac{n}{\log p}$. Furthermore, it is shown that this is the only regime where adaptive confidence intervals exist, even over two given parameter spaces. For example, when $k_1 \ll \frac{\sqrt{n}}{\log p}$ and $k_1\ll k_2$, it is impossible to construct  a confidence interval for the $\ell_2$ loss with guaranteed coverage probability over $\Theta_0(k_2)$ (consequently also over $\Theta_0(k_1)$) and with the expected length automatically adjusted to the sparsity. Similarly, for the $\ell_q$ loss with $1\leq q<2$,  adaptive confidence intervals is impossible over $\Theta_0(k_1)$ and $\Theta_0(k_2)$ for $k_1\ll k_2 \lesssim \frac{n}{\log p}$. Regarding confidence intervals for the $\ell_q$ loss with
$1\leq q\leq 2$ in the case of unknown $\Sigma$ and $\sigma$, the impossibility of adaptivity also holds over $\Theta(k_1)$ and $\Theta(k_2)$ for $k_1\ll k_2 \lesssim \frac{n}{\log p}$.

Establishing rate-optimal lower bounds 
requires the development of new technical tools. One main difference between loss estimation and the traditional parameter estimation is that for loss estimation  the constraint is on the performance of the estimator $\widehat{\beta}$ of the regression vector $\beta$, but the lower bound is on the difficulty of estimating its loss $\|\widehat{\beta} - \beta\|_q^2$.
We introduce useful new lower bound techniques for the minimax estimation error  and the expected length of adaptive confidence intervals for the  loss $\|\widehat{\beta} - \beta\|_q^2$. In several important cases, it is necessary to test a composite null against a composite alternative in order to establish rate sharp lower bounds. The technical tools developed in this paper can also be of independent interest.

In addition to $\Theta_0(k)$ and $\Theta(k)$, we also study an intermediate parameter space where the noise level $\sigma$ is known and the design covariance matrix $\Sigma$ is unknown but of certain structure. { Lower bounds for the expected length of minimax and adaptive confidence intervals for $\|\widehat{\beta}-\beta\|_q^2$ over this parameter space are established for a broad collection of estimators $\widehat{\beta}$ and are shown to be rate sharp for the class of rate-optimal estimators.} 
Furthermore, the lower bounds developed in this paper have wider implications. In particular, it is shown that they lead immediately to minimax lower bounds for estimating $\|\beta\|_q^2$ and the expected length of confidence intervals for $\|\beta\|_q^2$ with $1\le q\le 2$.



\subsection{Comparison with other works}

Statistical inference on the loss of specific estimators of $\beta$ has been considered in the recent literature. 
The papers \citep{donoho2011noise,bayati2012lasso} established, in the setting $\Sigma=\rm I$ and ${n}/{p} \rightarrow \delta\in (0,\infty)$,  the limit of the normalized loss $\frac{1}{p}\|\widehat{\beta}(\lambda)-\beta\|_2^2$ where $\widehat{\beta}(\lambda)$ is the Lasso estimator with a pre-specified tuning parameter $\lambda$. Although \citep{donoho2011noise,bayati2012lasso} provided an exact asymptotic expression of the normalized loss, the limit itself depends on the unknown $\beta$.  In a similar setting, the paper \citep{thrampoulidis2015asymptotically} established the limit of a normalized $\ell_2$ loss of the square-root Lasso estimator.
These limits of the normalized losses help understand the properties of the corresponding estimators of $\beta$, but they do not lead to an estimate of the loss.
Our results imply that although these normalized losses have a limit under some regularity conditions, such losses cannot be estimated well in most settings.

A recent paper, \citep{janson2015eigenprism}, constructed a confidence interval for $\|\widehat{\beta}-\beta\|_2^2$ in the case of known $\Sigma={\rm I}$, unknown noise level $\sigma$, and moderate dimension where $n/p\rightarrow \xi \in (0, 1)$ and no sparsity is assumed on $\beta$. 
While no sparsity assumption on $\beta$ is imposed, their method requires the assumption of $\Sigma={\rm I}$ and  $n/p \rightarrow \xi \in (0,1)$.  
In contrast, in this paper, we consider both unknown $\Sigma$ and known $\Sigma={\rm I}$ settings, while allowing $p\gg n$ and assuming sparse $\beta$. 


Honest adaptive inference has been studied in the nonparametric function estimation literature, including \citep{cai2004adaptation}  for adaptive confidence intervals for linear functionals, \citep{hoffmann2011adaptive,cai2014adaptive} for adaptive confidence bands, and \citep{Cai06,robins2006adaptive} for adaptive confidence balls, and in the high-dimensional linear regression literature, including \citep{nickl2013confidence} for adaptive confidence set  and \citep{cai2015regci} for adaptive confidence interval for linear functionals. 
In this paper, we develop new lower bound tools, Theorems \ref{thm: lower bound loss estimation generalization} and \ref{thm: lower bound for CI of risk high prob known variance}, to establish the possibility of adaptive confidence intervals for $\|\widehat{\beta}-\beta\|_q^2$. The connection between $\ell_2$ loss considered in the current paper and the work \citep{nickl2013confidence} is discussed in more detail in Section \ref{sec: l2 loss}. 
 
\subsection{Organization}
Section \ref{sec: minimax global loss} establishes the minimax lower bounds of estimating the loss $\|\widehat{\beta}-\beta\|_{q}^2$ with $1\leq q\leq 2$ over both $\Theta_0(k)$ and $\Theta(k)$ and shows that these bounds are rate sharp for the Lasso and scaled Lasso estimators, respectively. We then turn to interval estimation of  $\glosss$. Sections \ref{sec: CI for known} and  \ref{sec: CI for unknown} present the minimax and adaptive minimax lower bounds for the expected length of confidence intervals for $\glosss$ over $\Theta_0(k)$ and $\Theta(k)$. For Lasso and scaled Lasso estimators, we show that the lower bounds can be achieved and investigate the possibility of adaptivity. Section \ref{sec: rate optimal estimators} considers the rate-optimal estimators and establishes the minimax convergence rate of estimating their $\ell_{q}$ losses. Section \ref{sec: lower bound tool} presents new minimax lower bound techniques for estimating the loss $\glosss$. 
Section \ref{sec: more general parameter space} discusses the minimaxity and adaptivity in another setting, where the noise level $\sigma$ is known and the design covariance matrix $\Sigma$ is unknown but of certain structure. Section \ref{sec: lq norm functional} applies the newly developed lower bounds to establish lower bounds for a related problem, that of estimating $\|\beta\|_q^2$. Section \ref{sec: proof} proves the main results and additional proofs are given in the supplemental material \citep{cai2016asupplement}.

\subsection{Notation} 
\label{sec:notation}

For a matrix $X\in \R^{n\times p}$, $X_{i\cdot}$,  $X_{\cdot j}$, and $X_{i,j}$ denote respectively the $i$-th row,  $j$-th column, and  $(i,j)$ entry of the matrix $X$. 
 For a subset $J\subset\{1,2,\cdots,p\}$, $|J|$ denotes the cardinality of $J$, $J^{c}$ denotes the complement $\{1,2,\cdots,p\} \backslash J$, $X_{J}$ denotes the submatrix of $X$ consisting of columns  $X_{\cdot j}$ with $j\in J$ and for a vector $x\in \R^{p}$, $x_{J}$ is the subvector of $x$ with indices in $J$.
For a vector $x\in \R^{p}$, $\suppx$ denotes the support of $x$ and
the $\ell_q$ norm of $x$ is defined as $\|x\|_{q}=\left(\sum_{i=1}^{p}|x_i|^q\right)^{\frac{1}{q}}$ for $q \geq 0$ with $\|x\|_0=|\suppx|$ and $\|x\|_{\infty}=\max_{1\leq j \leq p}|x_j|$. 
For $a\in \R$, $a_{+}=\max\left\{a,0\right\}$. We use $\max\|\Xj\|_2$ as a shorthand for $\max_{1\leq j \leq p}\|\Xj\|_2$ and $\min\|\Xj\|_2$ as a shorthand for $\min_{1\leq j \leq p}\|\Xj\|_2$. For a matrix $A$, we define the spectral norm $\|A\|_{2}=\sup_{\|x\|_2 = 1} \|Ax\|_2$ and the matrix $\ell_1$ norm $\|A\|_{L_1}=\sup_{1\leq j\leq p}\sum_{i=1}^{p}|A_{ij}|$; For a symmetric matrix $A$, $\lambda_{\min}\left(A\right)$ and $\lambda_{\max}\left(A\right)$  denote respectively the smallest and largest eigenvalue of $A$. We use $c$ and $C$ to denote generic positive constants that may vary from place to place. For two positive sequences $a_n$ and $b_n$,  $a_n \lesssim b_n$ means $a_n \leq C b_n$ for all $n$ and $a_n \gtrsim b_n $ if $b_n\lesssim  a_n$ and $a_n \asymp b_n $ if $a_n \lesssim b_n$ and $b_n \lesssim a_n$, and $a_n \ll b_n$ if $\limsup_{n\rightarrow\infty} \frac{a_n}{b_n}=0$ and $a_n \gg b_n$ if $b_n \ll a_n$. 

\section{Minimax estimation of the $\ell_{q}$ loss}
\label{sec: minimax global loss}

We begin by presenting the minimax framework for estimating the $\ell_{q}$ loss, $\glosss$, of a given estimator $\widehat{\beta}$, and then establish the minimax lower bounds for the estimation error  for a broad collection of estimators $\widehat{\beta}$. 
We also show that such minimax lower bounds can be achieved for the Lasso and scaled Lasso estimators. 

\subsection{Problem formulation}
Recall the high-dimensional linear model,
\begin{equation}
y_{n\times 1}=X_{n\times p}\beta_{p\times 1}+\epsilon_{n\times 1}, \quad \epsilon \sim N_n(0,\sigma^2 {\rm I}).
\label{eq: linear model}
\end{equation}
We focus on the random design with $X_{i\cdot} \stackrel{iid}{\sim} N\left(0,\Sigma\right)$ and $X_{i\cdot}$ and $\epsilon_i$ are independent. 
Let $Z=(X,y)$ denote the observed data and $\widehat \beta$ be a given estimator of $\beta$. Denoting by $\lest(Z)$ any estimator of the loss $\glosss$, the minimax rate of convergence for estimating $\glosss$ over a parameter space $\Theta$ is defined as the largest quantity $\gamma_{\widehat{\beta},\ell_{q}}(\Theta)$ such that
\begin{equation}
\inf_{\lest} \sup_{\theta \in \Theta} \PP_{\theta}\left(|\lest(Z)-\glosss| \geq \gamma_{\widehat{\beta},\ell_{q}}(\Theta)\right) \geq \delta, 
\label{eq: minimax risk estimation for specific estimator high prob}
\end{equation}
 for some constant $\delta>0$ not depending on $n$ or $p$. We shall write $\lest$ for $\lest(Z)$ when there is no confusion.



We denote the parameter by $\theta=\left(\beta,\Sigma,\sigma\right)$, which consists of the signal $\beta$, the design covariance matrix $\Sigma$ and the noise level $\sigma$. For a given $\theta=\left(\beta,\Sigma,\sigma\right)$, we use $\beta(\theta)$ to denote the corresponding $\beta$.
Two settings are considered: The first is known design covariance matrix $\Sigma={\rm I}$ and known noise level $\sigma=\sigma_0$ and the other is unknown $\Sigma$ and $\sigma$. In the first setting, we consider  the following parameter space that  consists of $k$-sparse signals,
\begin{equation}
\Theta_0(k)=\left\{\left(\beta,{\rm I},\sigma_0 \right): \|\beta\|_0\leq k\right\}, 
\label{eq: parameter space with fixed variance and known design}
\end{equation}
and in the second setting, we consider
\begin{equation}
\Theta(k)=\left\{\left(\beta,\Sigma,\sigma\right): \|\beta\|_0\leq k, \; \frac{1}{\Md}\leq \lambda_{\min} \left(\Sigma\right)\leq \lambda_{\max} \left(\Sigma\right) \leq \Md, \; 0<\sigma\leq \Mv \right\},
\label{eq: parameter space with unknown variance}
\end{equation}
where $\Md\geq 1$ and $\Mv>0$ are constants. The parameter space $\Theta_0(k)$ 
is a subset of $\Theta(k)$, which consists of $k$-sparse signals with unknown $\Sigma$ and $\sigma$.  

The minimax rate $\gamma_{\widehat{\beta},\ell_{q}}(\Theta)$ for estimating $\glosss$ also 
depends on the estimator $\widehat{\beta}$. Different estimators $\widehat{\beta}$ could lead to different losses $\glosss$ and in general the difficulty of estimating the loss $\glosss$ varies with $\widehat{\beta}$.
We first recall the properties of some state-of-art estimators and then specify the collection of estimators on which we focus in this paper.
As shown in \citep{candes2007dantzig,bickel2009simultaneous,belloni2011square,sun2012scaled}, Lasso, Dantzig Selector, scaled Lasso and square-root Lasso satisfy the following property if the tuning parameter is properly chosen,
\begin{equation}
\sup_{\theta \in \Theta\left(k\right)} \PP_{\theta}\left(\glosss \geq C k^{\frac{2}{q}} \frac{\log p}{n} 
\right)\rightarrow 0,
\label{eq: rate optimal adaptive loss}
\end{equation}
where $C>0$ is a constant. The minimax lower bounds  established in \citep{verzelen2012minimax, raskutti2011minimax,ye2010rate} imply that $k^{\frac{2}{q}} \frac{\log p}{n} $ is the optimal rate for estimating $\beta$ over the parameter space $\Theta(k)$. It should be stressed that all of these algorithms do not require knowledge of the sparsity $k$ and  is thus adaptive to the sparsity provided $k \lesssim \frac{n}{\log p}$.
We consider a broad collection of estimators $\widehat{\beta}$ satisfying one of the following  two assumptions.
\begin{enumerate}
\item[(A1)] The estimator $\widehat{\beta}$ satisfies, for some $\theta_0=\left(\betanull,{\rm I},\sigma_0 \right)\in \Theta_0(k)$, 
\begin{equation}
\PP_{\theta_0} \left(\|\widehat{\beta}-\betanull\|_{q}^2\geq C^{*} \|\betanull\|_0^{\frac{2}{q}} {\frac{\log p}{n}}\sigma^2_0\right) \leq \alpha_0,
 \label{eq: gloss adaptivity assumption}
 \end{equation}
 where $0\leq \alpha_0<\frac{1}{4}$ and $C^{*}>0$ are constants.
\item[(A2)] The estimator  $\widehat{\beta}$ satisfies 
\begin{equation}
{ \sup_{\left\{\theta=\left(\betanull ,{\rm I}, \sigma\right):\sigma\leq 2\sigma_0\right\}}}\PP_{\theta}\left(\|\widehat{\beta}-\betanull\|_q^2\geq C^{*} \|\betanull\|_0^{\frac{2}{q}} {\frac{\log p}{n}} \sigma^2 \right)\leq \alpha_0,
\label{eq: adaptation high probability lq difference}
\end{equation}
where $0\leq \alpha_0<\frac{1}{4}$ and $C^{*}>0$ are constants and $\sigma_0>0$ is given.
\end{enumerate}

In view of the minimax rate given in \eqref{eq: rate optimal adaptive loss}, Assumption $({\rm A}1)$ requires $\widehat{\beta}$ to be a good estimator of $\beta$ at at least one point $\theta_0\in \Theta_0(k)$. Assumption $({\rm A}2)$ is slightly stronger than $({\rm A}1)$ and requires $\widehat{\beta}$ to estimate $\beta$ well for a single $\beta^*$ but over a range of noise levels $\sigma\le 2 \sigma_0$ while $\Sigma=\rm I$. Of course,  any estimator $\widehat{\beta}$ satisfying \eqref{eq: rate optimal adaptive loss} satisfies both $({\rm A}1)$  and $({\rm A}2)$. 
 In addition to Assumptions (A1) and (A2), we also introduce the following sparsity assumptions that will be used in various theorems.
\begin{enumerate}
\item[(B1)] Let $c_0$ be the constant defined in \eqref{eq: key constant}. The sparsity levels $k$ and $k_0$ satisfy $k \leq c_0 \min\{p^{\gamma}, \frac{n}{\log p}\}$ for some constant $0\leq \gamma <\frac{1}{2}$ and $k_0\leq c_0 \min\{k, \frac{\sqrt{n}}{\log p}\}$.  

\item[(B2)] 
The sparsity levels $k_1,k_2$ and $k_0$ satisfy $k_1\leq k_2 \leq c_0 \min\{p^{\gamma}, \frac{n}{\log p}\}$ for  some constant $0\leq \gamma <\frac{1}{2}$ and $c_0>0$ and $k_0 \leq c_0\min\{k_1,\frac{\sqrt{n}}{\log p}\}$. 
\end{enumerate}

\subsection{Minimax estimation of the $\ell_{q}$ loss over $\Theta_0(k)$}

The following theorem establishes the minimax lower bounds for estimating the  loss  $\glosss$ over the parameter space $\Theta_{0}\left(k\right)$.

\begin{Theorem}
\label{thm: estimating the global loss known}
Suppose that the sparsity levels $k$ and $k_0$ satisfy Assumption $({\rm B}1)$. For any estimator $\widehat{\beta}$ satisfying Assumption $({\rm A}1)$ with $\|\betanull\|_0\leq k_0$, 
\begin{equation}
\inf_{\widehat{L}_2} \sup_{\theta \in \Theta_0(k)} \PP_{\theta}\left(|\widehat{L}_2 - \|\widehat{\beta}-\beta\|_2^2| \geq c \min\left\{k{\frac{\log p}{n}},\frac{1}{\sqrt{n}}\right\}\sigma_0^2\right) \geq \delta.
\label{eq: key lower bound in high prob known application q=2}
\end{equation}
For any estimator $\widehat{\beta}$ satisfying Assumption $({\rm A}2)$ with $\|\betanull\|_0\leq k_0$, 
\begin{equation}
\inf_{\lest} \sup_{\theta \in \Theta_0(k)} \PP_{\theta}\left(|\lest-\glosss| \geq c k^{\frac{2}{q}}{\frac{\log p}{n}}\sigma_0^2\right) \geq \delta, \quad \text{for} \; 1\leq q<2,
\label{eq: key lower bound in high prob known application general q}
\end{equation}
where $\delta>0$ and $c>0$ are constants.
\end{Theorem}
\begin{Remark}  {\rm 
Assumption  (A1)  restricts our focus to estimators that can perform well at at least one point $\left(\betanull,{\rm I},\sigma_0\right) \in \Theta_0(k)$. This weak condition makes the established lower bounds widely applicable as the benchmark for evaluating estimators of the $\ell_q$ loss of any $\widehat{\beta}$ that performs well at a proper subset, or even a single point of the whole parameter space. 

In this paper, we focus on estimating the loss $\|\widehat{\beta}-\beta\|_q^2$ with $1\leq q\leq 2$. Similar results can be established for the loss in the form of $\|\widehat{\beta}-\beta\|_q^q$ with $1\leq q\leq 2$; Under the same assumptions as those in Theorem \ref{thm: estimating the global loss known}, the lower bounds for estimating the loss $\|\widehat{\beta}-\beta\|_q^q$ hold with replacing the convergence rates with their $\frac{q}{2}$ power; that is, \eqref{eq: key lower bound in high prob known application q=2} remains the same while the convergence rate $k^{\frac{2}{q}}(\sqrt{{\log p}/{n}}\sigma_0)^2$ in \eqref{eq: key lower bound in high prob known application general q} is replaced by $k(\sqrt{{\log p}/{n}}\sigma_0)^{q}$. Similarly,  all the results established in the rest of the paper for $\|\widehat{\beta}-\beta\|_q^2$ hold for $\|\widehat{\beta}-\beta\|_q^q$ with corresponding convergence rates replaced by their $\frac{q}{2}$ power.
}
\end{Remark}
Theorem \ref{thm: estimating the global loss known} establishes the minimax lower bounds for estimating the $\ell_2$ loss $\|\widehat{\beta}-\beta\|_2^2$ of any estimator $\widehat{\beta}$ satisfying Assumption $({\rm A}1)$ and the $\ell_q$ loss  $\|\widehat{\beta}-\beta\|_q^2$ with $1\le q<2$ of any estimator $\widehat{\beta}$ satisfying Assumption $({\rm A}2)$. We will take the Lasso estimator as an example and demonstrate the implications of the above theorem. 
We randomly split $Z=(y,X)$ into subsamples $Z^{\left(1\right)}=\left(y^{\left(1\right)},X^{\left(1\right)}\right)$ and $Z^{\left(2\right)}=\left(y^{\left(2\right)},X^{\left(2\right)}\right)$ with sample sizes $n_1$ and $n_2$, respectively. The Lasso estimator $\betal$ based on the first subsample $Z^{\left(1\right)}=\left(y^{\left(1\right)},X^{\left(1\right)}\right)$ is defined as
\begin{equation}
\betal=\argmin_{\beta \in \R^{p}}\frac{\|\ya-\xa\beta\|_2^2}{n_1}+\lambda \sum_{j=1}^{p} \frac{\|\xa_{\cdot j}\|_2}{\sqrt{n_1}} |\beta_j|,
\label{eq: Lasso}
\end{equation}
where $\lambda=A\sqrt{{\log p}/{\nf}}\sigma_0$ with $A>\sqrt{2}$ being a pre-specified constant. Without loss of generality, we assume $\nf\asymp \ns.$ 
For the case $1\leq q <2$, \eqref{eq: rate optimal adaptive loss} and \eqref{eq: key lower bound in high prob known application general q}  together imply that the estimation of the $\ell_{q}$ loss $\|\betal-\beta\|_{q}^2$ is impossible since the lower bound can be achieved by the trivial estimator of the loss, 0. That is,
$\sup_{\theta \in \Theta_0(k)} \PP_{\theta}\left(|0-\|\betal-\beta\|_{q}^2| \geq C k^{\frac{2}{q}} {\frac{\log p}{n}}\right) \rightarrow 0.$

For the case $q=2$, in the regime $k\ll\frac{\sqrt{n}}{\log p}$, the lower bound $\frac{k \log p}{n}$ in \eqref{eq: key lower bound in high prob known application q=2} can be achieved by the zero estimator and hence estimation of the loss $\|\betal-\beta\|_2^2$ is impossible. However, the interesting case is when $\frac{\sqrt{n}}{\log p} \lesssim k\lesssim \frac{n}{\log p}$,
the loss estimator $\widetilde{L}_2$ proposed in \eqref{eq: proposed estimator of l2 loss estimation} achieves the minimax lower bound $\frac{1}{\sqrt{n}}$ in \eqref{eq: key lower bound in high prob known application q=2}, which cannot be achieved by the zero estimator. We now detail the construction of the loss estimator $\widetilde{L}_2$.
Based on the second half sample $Z^{\left(2\right)}=\left(y^{\left(2\right)},X^{\left(2\right)}\right)$, we propose the following estimator,
\begin{equation}
\widetilde{L}_2 =\left(\frac{1}{\ns}\left\|y^{\left(2\right)}-X^{\left(2\right)}\betal \right\|_2^2-\sigma_0^2\right)_{+}.
\label{eq: proposed estimator of l2 loss estimation}
\end{equation}
Note that the first subsample $Z^{\left(1\right)}=\left(y^{\left(1\right)},X^{\left(1\right)}\right)$  is used to produce the Lasso estimator $\betal$ in \eqref{eq: Lasso} and the second subsample $Z^{\left(2\right)}=\left(y^{\left(2\right)},X^{\left(2\right)}\right)$ is retained to evaluate the loss $\|\betal-\beta\|_2^2$. 
Such sample splitting technique is similar to cross-validation and has been used in \citep{nickl2013confidence} for constructing confidence sets for $\beta$ and in \citep{janson2015eigenprism} for confidence intervals for the $\ell_2$ loss. 

The following proposition establishes that the estimator $\widetilde{L}_2$ achieves the minimax lower  bound of  \eqref{eq: key lower bound in high prob known application q=2} over the regime $\frac{\sqrt{n}}{\log p} \lesssim k \lesssim \frac{n}{\log p}$.

\begin{Proposition}
\label{prop: positive estimation q=2}
Suppose that $k\lesssim \frac{n}{\log p}$ and $\betal$ is the Lasso estimator defined in \eqref{eq: Lasso} with $A>\sqrt{2}$, then the estimator of loss proposed in \eqref{eq: proposed estimator of l2 loss estimation} satisfies, for any sequence $\delta_{n,p}\rightarrow \infty$,
\begin{equation}
\limsupnp\sup_{\theta \in \Theta_0(k)} \PP_{\theta}\left(\left| \widetilde{L}_2 -\|\betal-\beta\|_2^2\right| \geq \delta_{n,p}\frac{1}{\sqrt{n}} \right) = 0.
\label{eq: upper bound known q=2}
\end{equation}
\end{Proposition}

\subsection{Minimax estimation of the $\ell_{q}$ loss over $\Theta(k)$}

We now turn to the case of unknown $\Sigma$ and $\sigma$ and establish the minimax lower bound for estimating the $\ell_{q}$ loss over the parameter space $\Theta(k)$.
\begin{Theorem}
\label{thm: estimating the global loss}
Suppose that the sparsity levels $k$ and $k_0$ satisfy Assumption $({\rm B}1)$. For any estimator $\widehat{\beta}$ satisfying Assumption $({\rm A}1)$
 with $\|\betanull\|_0\leq k_0$,  
 \begin{equation}
\inf_{\lest} \sup_{\theta \in \Theta(k)} \PP_{\theta}\left(|\lest-\glosss| \geq c k^{\frac{2}{q}} {\frac{\log p}{n}}\right) \geq \delta, \quad 1\leq q\leq 2,
\label{eq: key lower bound in high prob unknown application}
\end{equation}
where $\delta>0$ and $c>0$ are constants.
\end{Theorem}

Theorem \ref{thm: estimating the global loss} provides a minimax lower bound for estimating the $\ell_{q}$ loss of any estimator $\widehat{\beta}$ satisfying Assumption $({\rm A}1)$, including the scaled Lasso estimator defined as
\begin{equation}
\{\betasl,\hat{\sigma}\}=\argmin_{\beta \in \R^{p},\sigma
\in \R^{+}}\frac{\|y-X\beta\|_2^2}{2n\sigma}+\frac{\sigma}{2}+\lambda_0 \sum_{j=1}^{p} \frac{\|\Xj\|_2}{\sqrt{n}} |\beta_j|,
\label{eq: scaled Lasso}
\end{equation}
where $\lambda_0= A \sqrt{{\log p}/{n}}$ with $A>\sqrt{2}$.
%
Note that for the scaled Lasso estimator, the lower bound in \eqref{eq: key lower bound in high prob unknown application} can be achieved by the trivial loss estimator 0 in the sense,
$\sup_{\theta \in \Theta(k)} \PP_{\theta}\left(|0-\|\betasl-\beta\|_{q}^2| \geq C k^{\frac{2}{q}} {\frac{\log p}{n}}\right) \rightarrow 0,$
and hence estimation of loss is impossible in this case.  
\section{Minimaxity and adaptivity of confidence intervals over $\Theta_0(k)$}
\label{sec: CI for known}
We focused  in the last section on point estimation of the $\ell_q$ loss and showed the impossibility of loss estimation except for one  regime. The results naturally lead to another question: Is it possible to construct ``useful" confidence intervals for $\glosss$ that can provide non-trivial upper and lower bounds for the loss? In this section, after introducing the framework for minimaxity and adaptivity of confidence intervals, we consider the case of known $\Sigma={\rm I}$ and $\sigma=\sigma_0$ and establish the minimaxity and adaptivity lower bounds for the expected length of confidence intervals for the $\ell_q$ loss of a broad collection of estimators over the parameter space $\Theta_0(k)$. We also show that such minimax lower bounds can be achieved for the Lasso estimator and then discuss the possibility of adaptivity using the Lasso estimator as an example. 
The case of unknown $\Sigma$ and $\sigma$ will be the focus of  the next section.

\subsection{Framework for minimaxity and adaptivity of confidence intervals}
\label{framework.sec}
\begin{figure}[h!]
\centering
\includegraphics[width=3.4in,height=1.2in]{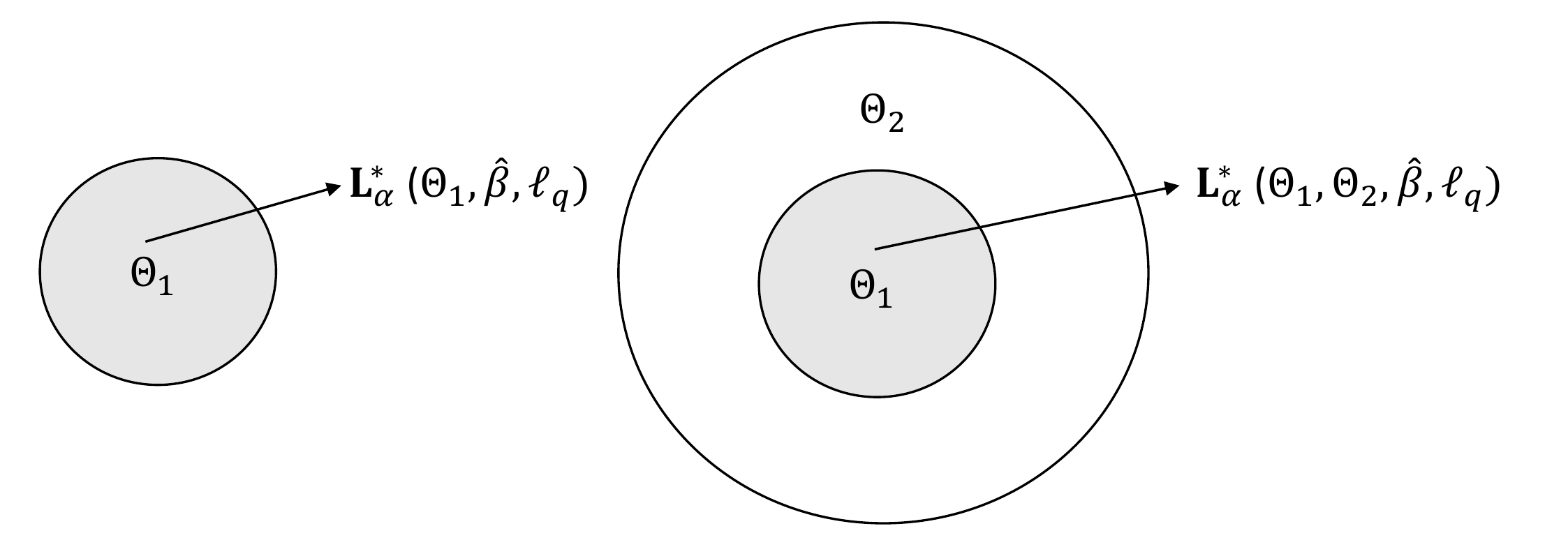}
\caption{The plot demonstrates the definition of $\RR_{\alpha}^{*}\left(\Thetas,\widehat{\beta},\ell_{q}\right)$ and $\RR_{\alpha}^{*}\left(\Thetas,\Thetal,\widehat{\beta},\ell_{q}\right)$.} 
\label{fig:def of adaptation parameter}
\end{figure}

In this section, we introduce the following decision theoretical framework for confidence intervals of the loss $\glosss$.  
Given $0<\alpha<1$ and the parameter space $\Theta$ and the loss $\glosss$, denote by  $\II_{\alpha}\left(\Theta, \widehat{\beta},\ell_{q}\right)$ the set of all $(1-\alpha)$ level confidence intervals for $\glosss$ over $\Theta$,
\begin{equation}
\II_{\alpha}\left(\Theta, \widehat{\beta},\ell_{q}\right)=\left\{\gCIZR=\left[l\left(Z\right),u\left(Z\right)\right]: \inf_{\theta \in \Theta} \PP_{\theta}\left(\gloss \in \gCIZR\right)\geq 1-\alpha \right\}.
\label{eq: def of set of CI}
\end{equation}
We will write ${\rm CI}_{\alpha}$ for $\gCIZR$ when there is no confusion.
For any confidence interval $\gCIZR=\left[l\left(Z\right),u\left(Z\right)\right]$, its length is denoted by $\RR\left(\gCIZR\right)=u\left(Z\right)-l\left(Z\right)$ and the maximum expected length over a parameter space $\Thetas$ is defined as
\begin{equation}
\RR\left(\gCIZR,\Thetas\right)=\sup_{\theta\in \Thetas} \E_{\theta} \RR\left(\gCIZR\right).
\label{eq: length over parameter space}
\end{equation}
 For two nested parameter spaces $\Thetas \subseteq \Thetal$, we define the benchmark 
 $\RR_{\alpha}^{*}\left(\Thetas,\Thetal,\widehat{\beta}, \ell_{q}\right)$, measuring the degree of adaptivity over the nested spaces $\Thetas \subset \Thetal$,
 \begin{equation}
\RR_{\alpha}^{*}\left(\Thetas,\Thetal,\widehat{\beta}, \ell_{q}\right)=\inf_{\gCIZR\in \II_{\alpha}\left(\Thetal,\widehat{\beta}, \ell_{q}\right)}\sup_{\theta \in \Thetas}\E_{\theta} \RR\left(\gCIZR\right).
\label{eq: benchmark loss}
\end{equation}
We will write $\RR_{\alpha}^{*}\left(\Thetas,\widehat{\beta},\ell_{q}\right)$ for $\RR_{\alpha}^{*}\left(\Thetas,\Thetas,\widehat{\beta},\ell_{q}\right)$, which is the minimax expected length of confidence intervals of $\glosss$ over $\Thetas$. The benchmark $\RR_{\alpha}^{*}\left(\Thetas,\Thetal,\widehat{\beta}, \ell_{q}\right)$ is the infimum of the maximum expected length over $\Thetas$ among all $(1-\alpha)$-level confidence intervals over $\Thetal$. In contrast, $\RR_{\alpha}^{*}\left(\Thetas,\widehat{\beta}, \ell_{q}\right)$ is considering  all $(1-\alpha)$-level confidence intervals over $\Thetas$. In words, if there is prior information that the parameter lies in the smaller parameter space $\Thetas$, $\RR_{\alpha}^{*}\left(\Thetas,\widehat{\beta}, \ell_{q}\right)$ measures the benchmark length of confidence intervals over the parameter space $\Thetas$, which is illustrated in the left of Figure \ref{fig:def of adaptation parameter}; however, if there is only prior information that the parameter lies in the larger parameter space $\Thetal$, $\RR_{\alpha}^{*}\left(\Thetas,\Thetal,\widehat{\beta}, \ell_{q}\right)$ measures the benchmark length of confidence intervals over the parameter space $\Thetas$, which is illustrated in the right of Figure \ref{fig:def of adaptation parameter}. 

Rigorously, we define a confidence interval ${\rm CI}^{*}$ to be simultaneously adaptive over $\Thetas$ and $\Thetal$ if ${\rm CI}^{*}\in \II_{\alpha}\left(\Thetal,\widehat{\beta},\ell_{q}\right),$
 \begin{equation}
\RR\left({\rm CI}^{*},\Thetas\right)\asymp \RR_{\alpha}^{*}\left(\Thetas,\widehat{\beta},\ell_{q}\right), \; \mbox{ and }\; \RR\left({\rm CI}^{*},\Thetal\right)\asymp \RR_{\alpha}^{*}\left(\Thetal,\widehat{\beta},\ell_{q}\right).
\label{eq: def of adaptivity for global loss}
\end{equation}
The condition \eqref{eq: def of adaptivity for global loss} means that the confidence interval ${\rm CI}^{*}$ has coverage over the larger parameter space $\Thetal$ and achieves the minimax rate over both $\Thetas$ and $\Thetal$.
Note that $\RR\left({\rm CI}^{*},\Thetas\right)\geq\RR_{\alpha}^{*}\left(\Thetas,\Thetal,\widehat{\beta}, \ell_{q}\right)$. If $\RR_{\alpha}^{*}\left(\Thetas,\Thetal,\widehat{\beta}, \ell_{q}\right)\gg \RR_{\alpha}^{*}\left(\Thetas,\widehat{\beta}, \ell_{q}\right),$ then the rate-optimal adaptation \eqref{eq: def of adaptivity for global loss} is impossible to achieve for $\Thetas \subset \Thetal$. 
Otherwise, it is possible to construct confidence intervals simultaneously adaptive over parameter spaces $\Thetas$ and $\Thetal$. The possibility of adaptation over parameter spaces $\Thetas$ and $\Thetal$ can thus be answered by
investigating the benchmark quantities 
 $\RR_{\alpha}^{*}\left(\Thetas,\widehat{\beta},\ell_{q}\right)$ and $\RR_{\alpha}^{*}\left(\Thetas, \Thetal,\widehat{\beta},\ell_{q}\right)$.
 Such framework has already been introduced in \citep{cai2015regci}, which studies the minimaxity and adaptivity of confidence intervals for linear functionals in high-dimensional linear regression.

We will adopt the minimax and adaptation framework discussed above and establish the minimax expected length $\RR_{\alpha}^{*}\left(\Theta_0(k), \widehat{\beta},\ell_{q}\right)$ and the adaptation benchmark $\RR_{\alpha}^{*}\left(\Theta_0(k_1), \Theta_0(k_2),\widehat{\beta},\ell_{q}\right)$. In terms of the minimax expected length and the adaptivity behavior, there exist fundamental differences between the case $q=2$ and $1\leq q<2$ . We will discuss them separately in the following two sections.

\subsection{Confidence intervals for the $\ell_2$ loss over $\Theta_0(k)$} 
\label{sec: l2 loss}
The following theorem establishes the minimax lower bound for the expected length of confidence intervals of $\|\widehat{\beta}-\beta\|_2^2$ over the parameter space  $\Theta_0(k)$.
\begin{Theorem}
\label{thm: minimax of CI for the l2 loss both known}
Suppose that $0<\alpha<\frac{1}{4}$ and the sparsity levels $k$ and $k_0$ satisfy Assumption $({\rm B}1)$. For any estimator $\widehat{\beta}$ satisfying Assumption $({\rm A}1)$ with $\|\betanull\|_0\leq k_0$,
then there is some constant $c>0$ such that 
\begin{equation}
\RR_{\alpha}^{*}\left(\Theta_0(k), \widehat{\beta},\ell_{2}\right) \geq  c \min\left\{\frac{k \log p}{n}, \frac{1}{\sqrt{n}}\right\}\sigma_0^2.
\label{eq: l2 lower bound}
\end{equation}
In particular, if $\betal$ is the Lasso estimator defined in \eqref{eq: Lasso} with $A>\sqrt{2}$, then the minimax expected length for $\left(1-\alpha\right)$ level confidence intervals of $\|\betal-\beta\|_2^2$ over $\Theta_0(k)$ is
\begin{equation}
\RR_{\alpha}^{*}\left(\Theta_0(k), \betal,\ell_{2}\right) \asymp \min\left\{\frac{k \log p}{n}, \frac{1}{\sqrt{n}}\right\}\sigma_0^2.
\label{eq: l2 minimax rate}
\end{equation}
\end{Theorem}

We now consider adaptivity of confidence intervals for the $\ell_2$ loss. The following theorem gives the lower bound for the benchmark $\RR_{\alpha}^{*}\left(\Theta_0(k_1),\Theta_0(k_2), \widehat{\beta},\ell_{2}\right)$.
We will then discuss Theorems \ref{thm: minimax of CI for the l2 loss both known} and \ref{thm: adaptation of CI for the l2 loss both known} together.





\begin{Theorem}
\label{thm: adaptation of CI for the l2 loss both known}
Suppose that $0<\alpha<\frac{1}{4}$ and the sparsity levels $k_1,k_2$ and $k_0$ satisfy Assumption $({\rm B}2)$. For any estimator $\widehat{\beta}$ satisfying Assumption $({\rm A}1)$ with $\|\betanull\|_0\leq k_0,$ then there is some constant $c>0$ such that
\begin{equation}
\RR_{\alpha}^{*}\left(\Theta_0(k_1),\Theta_0(k_2), \widehat{\beta},\ell_{2}\right) \geq c \min\left\{\frac{k_2\log p}{n} , \frac{1}{\sqrt{n}}\right\}\sigma_0^2.
\label{eq: l2 lower bound adaptation}
\end{equation}
In particular, if $\betal$ is the Lasso estimator defined in \eqref{eq: Lasso} with $A>\sqrt{2}$, the above lower bound can be achieved.
\end{Theorem}

The lower bound established in Theorem \ref{thm: adaptation of CI for the l2 loss both known}  implies that of Theorem \ref{thm: minimax of CI for the l2 loss both known} and  both lower bounds hold for a general class of estimators satisfying Assumption $({\rm A}1)$. There is a phase transition for the lower bound of the benchmark $\RR_{\alpha}^{*}\left(\Theta_0(k_1),\Theta_0(k_2), \widehat{\beta},\ell_{2}\right)$. In the regime $k_2\ll \frac{\sqrt{n}}{\log p}$, the lower bound in  \eqref{eq: l2 lower bound adaptation} is $\frac{k_2\log p}{n}\sigma_0^2$; when $\frac{\sqrt{n}}{\log p}\lesssim k_2\lesssim \frac{n}{\log p}$, the lower bound in \eqref{eq: l2 lower bound adaptation} is $\frac{1}{\sqrt{n}}\sigma_0^2$. For the Lasso estimator $\betal$ defined in \eqref{eq: Lasso}, the lower bound $\frac{k\log p}{n}\sigma_0^2$ in \eqref{eq: l2 lower bound} 
and $\frac{k_2 \log p}{n}\sigma_0^2$ in \eqref{eq: l2 lower bound adaptation} 
can be achieved by the confidence intervals ${\rm CI}_{\alpha}^{0} \left(Z,k,2\right)$ and ${\rm CI}_{\alpha}^{0} \left(Z,k_2,2\right)$ defined in \eqref{eq: theoretical CI for known case}, respectively. Such an interval estimator is also used for the $\ell_{q}$ loss with $1\leq q<2$. The minimax lower bound $\frac{1}{\sqrt{n}}\sigma_0^2$ in \eqref{eq: l2 minimax rate} and \eqref{eq: l2 lower bound adaptation} 
can be achieved by the following confidence interval,
\begin{equation}
{\rm CI}_{\alpha}^{1} \left(Z\right)= \left(\left(\frac{\psi\left(Z\right)}{\frac{1}{\ns}\chi^2_{1-\frac{\alpha}{2}}\left(\ns\right)}-\sigma_0^{2}\right)_{+},\left(\frac{\psi\left(Z\right)}{\frac{1}{\ns}\chi^2_{\frac{\alpha}{2}}\left(\ns\right)}-\sigma_0^{2}\right)_{+}\right), 
\label{eq: propose CI for square norm loss}
\end{equation}
where $\chi^2_{1-\frac{\alpha}{2}}\left(\ns\right)$ and $\chi^2_{\frac{\alpha}{2}}\left(\ns\right)$ are the $1-\frac{\alpha}{2}$ and $\frac{\alpha}{2}$ quantiles of $\chi^2$ random variable with $\ns$ degrees of freedom, respectively, and
\begin{equation}
\psi\left(Z\right)=\min\left\{\frac{1}{\ns}\left\|y^{\left(2\right)}-X^{\left(2\right)}\betal\right\|_2^2, \sigma_0^2\log p\right\}.
\label{eq: def of radius}
\end{equation}
Note that the two-sided confidence interval \eqref{eq: propose CI for square norm loss} is simply based on the observed data $Z$, not depending on any prior knowledge of the sparsity $k$. Furthermore, it is a two-sided confidence interval, which tells not only just an upper bound, but also a lower bound for the loss. The coverage property and the expected length of ${\rm CI}_{\alpha}^{1} \left(Z\right)$ are established in the following proposition.

\begin{Proposition}
\label{prop: l2 positive CI}
Suppose $k\lesssim \frac{n}{\log p}$ and $\betal$ is the estimator defined in \eqref{eq: Lasso} with $A>\sqrt{2}$.  
Then ${\rm CI}_{\alpha}^{1}\left(Z\right)$ defined in \eqref{eq: propose CI for square norm loss} satisfies, 
\begin{equation}
\liminfnp \inf_{\theta\in \Theta_0\left(k\right)}\PP\left(\|\betal-\beta\|_2^2\in {\rm CI}_{\alpha}^{1}\left(Z\right)\right)\geq 1-\alpha,
\label{eq: coverage of CI1 q=2}
\end{equation}
and
\begin{equation}
\RR\left({\rm CI}_{\alpha}^{1}\left(Z\right),\Theta_0\left(k\right)\right)\lesssim \frac{1}{\sqrt{n}}\sigma_0^2.
\label{eq: length of CI1 q=2}
\end{equation}
\end{Proposition}

\begin{figure}[h!]
\centering
\includegraphics[scale=0.38]{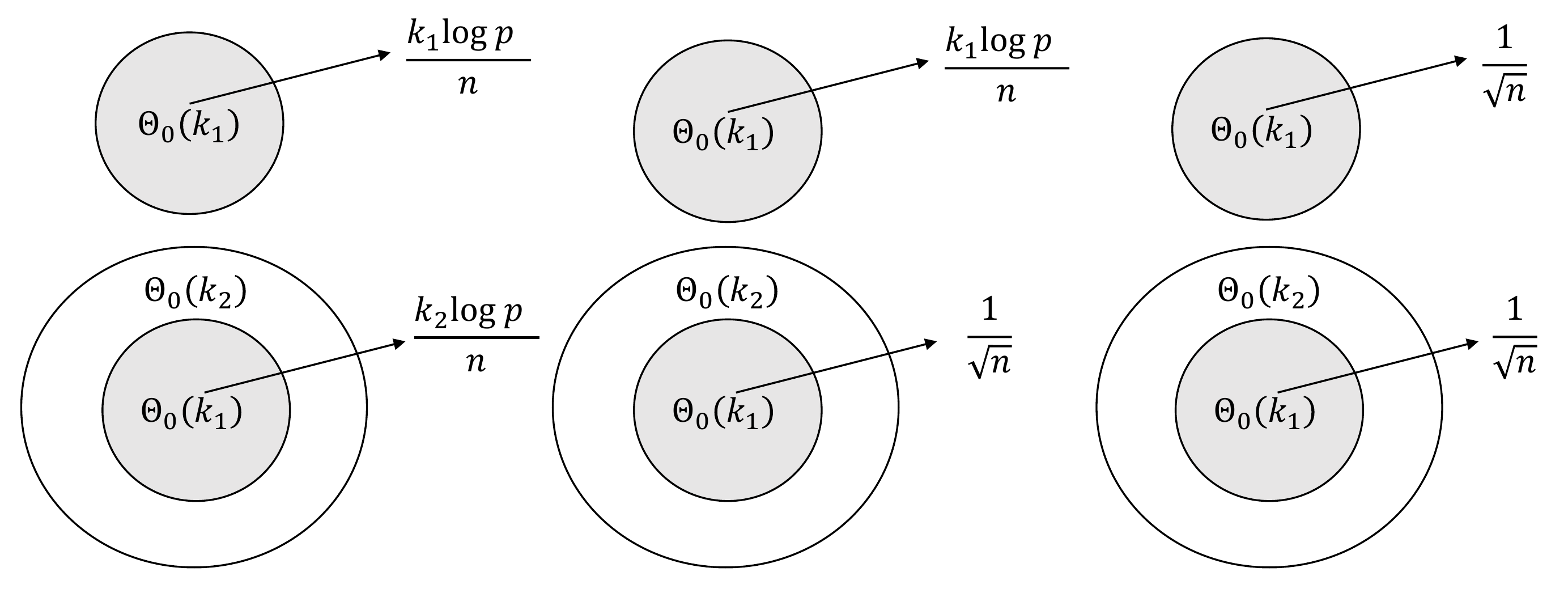}
\caption{Illustration of $\RR_{\alpha}^{*}\left(\Theta_0(k_1), \betal,\ell_{2}\right)$ (top) and   $\RR_{\alpha}^{*}\left(\Theta_0(k_1),\Theta_0(k_2), \betal,\ell_{2}\right)$(bottom) over regimes $\regimea$ (leftmost), $\regimeb$(middle) and $\regimec$ (rightmost).}
\label{fig:adaptivity behavior of l2 norm known}
\end{figure}

Regarding the Lasso estimator $\betal$ defined in \eqref{eq: Lasso}, we will discuss the possibility of adaptivity of confidence intervals for $\|\betal-\beta\|_2^2$.
The adaptivity behavior of confidence intervals for $\|\betal-\beta\|_2^2$ is demonstrated in Figure \ref{fig:adaptivity behavior of l2 norm known}. As illustrated in the rightmost plot of Figure \ref{fig:adaptivity behavior of l2 norm known}, in the regime $\regimec$, we obtain $\RR_{\alpha}^{*}\left(\Theta_0(k_2),\Theta_0(k_2), \betal,\ell_{2}\right)\asymp \RR_{\alpha}^{*}\left(\Theta_0(k_1), \betal,\ell_{2}\right) \asymp \frac{1}{\sqrt{n}}$, which implies that adaptation is possible over this regime. As shown in Proposition \ref{prop: l2 positive CI}, the confidence interval ${\rm CI}_{\alpha}^{1} \left(Z\right)$ defined in \eqref{eq: propose CI for square norm loss} is fully adaptive over the regime $\frac{\sqrt{n}}{\log p}\lesssim k\lesssim \frac{n}{\log p}$ in the sense of  \eqref{eq: def of adaptivity for global loss}.

Illustrated in the leftmost and middle plots of Figure \ref{fig:adaptivity behavior of l2 norm known},
it is impossible to construct an adaptive confidence interval for $\|\betal-\beta\|_2^2$ 
over regimes $\regimea$ and $k_1\ll\frac{\sqrt{n}}{\log p}\lesssim k_2\lesssim \frac{n}{\log p}$ since
$\RR_{\alpha}^{*}\left(\Theta_0(k_1),\Theta_0(k_2), \betal,\ell_{2}\right)\gg \RR_{\alpha}^{*}\left(\Theta_0(k_1), \betal,\ell_{2}\right)$ if $k_1\ll \frac{\sqrt{n}}{\log p}$ {and}  $k_1 \ll k_2.$
To sum up, adaptive confidence intervals for $\|\betal-\beta\|_2^2$ is only possible over the regime $\frac{\sqrt{n}}{\log p}\lesssim k \lesssim \frac{n}{\log p}$.

\subsubsection*{Comparison with confidence balls}

\begin{figure}[h!]
\centering
\includegraphics[width=3.6in,height=1in]{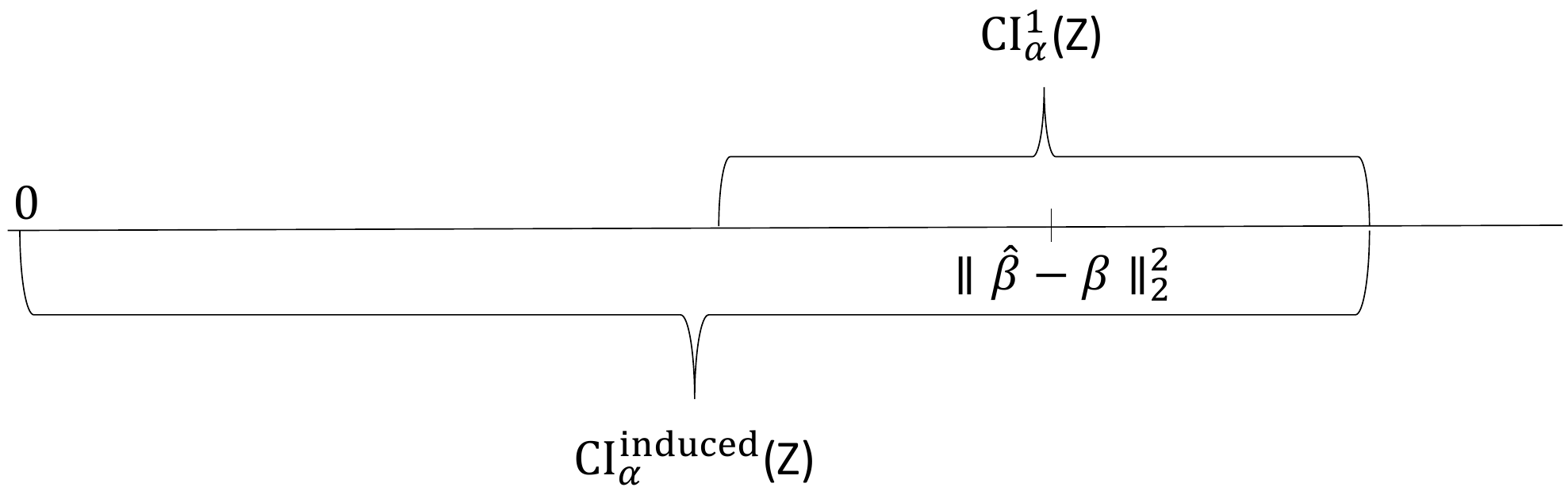}
\caption{Comparison of the two-sided confidence interval ${\rm CI}_{\alpha}^{1} \left(Z\right)$ with the one-sided confidence interval ${\rm CI}_{\alpha}^{\rm induced}\left(Z\right)$.}
\label{fig:comparisonofCI}
\end{figure}

We should note that the problem of constructing confidence intervals for $\|\widehat{\beta}-\beta\|_2^2$ is related to but different from that of constructing confidence sets for $\beta$ itself. Confidence balls constructed in \citep{nickl2013confidence} are of form $\left\{\beta: \; \|{\beta}-\widehat{\beta}\|_2^2 \leq u_{n}\left(Z\right)\right\}$, where $\widehat{\beta}$ can be the Lasso estimator and $u_{n}\left(Z\right)$ is a data dependent squared radius. See \citep{nickl2013confidence} for further details. A naive application of this confidence ball leads to a one-sided confidence interval for the loss $\|\widehat{\beta}-\beta\|_2^2$,
\begin{equation}
{\rm CI}_{\alpha}^{\rm induced}\left(Z\right)=\left\{\|\widehat{\beta}-{\beta}\|_2^2: \; \|\widehat{\beta}-{\beta}\|_2^2\leq u_{n}\left(Z\right)\right\}.
\label{eq: confidence ball introduced}
\end{equation}
Due to the reason that confidence sets for $\beta$ were sought for in Theorem 1 in \citep{nickl2013confidence},  confidence sets in the form $\left\{\beta: \; \|{\beta}-\widehat{\beta}\|_2^2 \leq u_{n}\left(Z\right)\right\}$ will suffice to achieve the optimal length. However, since our goal is to characterize $\|\widehat{\beta}-\beta\|_{2}^2$, we apply the unbiased risk estimation discussed in Theorem 1 of \citep{nickl2013confidence} and construct the two-sided confidence interval in \eqref{eq: propose CI for square norm loss}.
 Such a two-sided confidence interval is more informative than the one-sided confidence interval \eqref{eq: confidence ball introduced} since the one-sided confidence interval does not contain the information whether the loss is close to zero or not. Furthermore, as shown in \citep{nickl2013confidence}, the length of confidence interval ${\rm CI}_{\alpha}^{\rm induced}\left(Z\right)$ over the parameter space $\Theta_{0}(k)$ is of order  $\frac{1}{\sqrt{n}}+\frac{k \log p}{n}$.
The two-sided confidence interval ${\rm CI}_{\alpha}^{1} \left(Z\right)$ constructed in \eqref{eq: propose CI for square norm loss} is of expected length $\frac{1}{\sqrt{n}}$, which is much shorter than  $\frac{1}{\sqrt{n}}+\frac{k \log p}{n}$ in the regime $k \gg \frac{\sqrt{n}}{\log p}$. That is, the two-sided confidence interval \eqref{eq: propose CI for square norm loss} provides a more accurate interval estimator of the $\ell_2$ loss. This is illustrated in Figure \ref{fig:comparisonofCI}.


The lower bound technique developed in the literature of adaptive confidence sets \citep{nickl2013confidence} can also be used to establish some of the lower bound results for the case $q=2$ given in the present paper. However, new techniques are needed in order to establish the rate sharp lower bounds for the minimax estimation error \eqref{eq: key lower bound in high prob known application general q} in the region $\frac{\sqrt{n}}{\log p}\leq k\lesssim \frac{{n}}{\log p}$ and for the expected length of the confidence intervals \eqref{eq: lq lower bound adaptation} and \eqref{eq: l2 lower bound adaptation general} in the region $\regimec$, where it is necessary to test a composite null against a composite alternative in order to establish rate sharp lower bounds. 

\subsection{Confidence intervals for the $\ell_q$ loss with $1\leq q<2$ over $\Theta_0(k)$}
\label{sec: lq loss}
We now consider the case $1\leq q<2$ and  investigate the minimax expected length and adaptivity of confidence intervals for $\glosss$ over the parameter space $\Theta_0(k)$. 
The following theorem characterizes the minimax convergence rate for the expected length of confidence intervals. 
\begin{Theorem}
\label{thm: minimax of CI for the lq loss known variance}
Suppose that $0<\alpha<\frac{1}{4}$, $1\leq q< 2$ and the sparsity levels $k$ and $k_0$ satisfy Assumption $({\rm B}1)$. For any estimator $\widehat{\beta}$ satisfying Assumption $({\rm A}2)$ with $\|\betanull\|_0\leq k_0$, then there is some constant $c>0$ such that
\begin{equation}
\RR_{\alpha}^{*}\left(\Theta_0(k), \widehat{\beta},\ell_{q}\right) \geq  c k^{\frac{2}{q}}{\frac{\log p}{n}}\sigma_0^2.
\label{eq: lq lower bound}
\end{equation}
In particular, if $\betal$ is the Lasso estimator defined in \eqref{eq: Lasso} with $A>4\sqrt{2}$, 
then the minimax expected length for $(1-\alpha)$ level confidence intervals of $\|\betal-\beta\|_{q}^2$ over $\Theta_0(k)$ is
\begin{equation}
\RR_{\alpha}^{*}\left(\Theta_0(k), \betal,\ell_{q}\right) \asymp k^{\frac{2}{q}}{\frac{\log p}{n}}\sigma_0^2.
\label{eq: lq minimax}
\end{equation}
 \end{Theorem}
We now construct the confidence interval achieving the minimax convergence rate of \eqref{eq: lq minimax},
\begin{equation}
{\rm CI}_{\alpha}^{0} \left(Z,k,q\right)=\left(0,C^{*}(A,k) k^{\frac{2}{q}}\frac{\log p}{n} \right), 
\label{eq: theoretical CI for known case}
\end{equation}
where $C^{*}(A,k)=\max\left\{\frac{\left(22 A {\sigma}_0\right)^2}{\left(\frac{1}{4}-42\sqrt{\frac{2k \log p}{\nf}}\right)^4},\frac{\left(\frac{3\eta_0}{\eta_0+1} A \sigma_0\right)^2}{\left(\frac{1}{4}-(9+11\eta_0)\sqrt{\frac{2k \log p}{\nf}}\right)^4}\right\}$ with $\eta_0=1.01\frac{\sqrt{A}+\sqrt{2}}{\sqrt{A}-\sqrt{2}}$.
The following proposition establishes the coverage property and the expected length of ${\rm CI}^{0}_{\alpha}\left(Z,k,q\right)$. 

\begin{Proposition}
\label{prop: construction of CI known}
Suppose $k  \lesssim \frac{n}{\log p}$ and $\betal$ is the estimator defined in \eqref{eq: Lasso} with $A>4\sqrt{2}$. 
For $1\leq q \leq 2$, the confidence interval ${\rm CI}^{0}_{\alpha}\left(Z,k,q\right)$ defined in \eqref{eq: theoretical CI for known case} satisfies
\begin{equation}
\liminfnp \inf_{\theta\in \Theta_0\left(k\right)}\PP_{\theta}\left(\|\widehat{\beta}-\beta\|_q^2\in {\rm CI}_{\alpha}^{0}\left(Z,k,q\right)\right) =1, 
\label{eq: coverage of CI known}
\end{equation}
and
\begin{equation}
\RR\left({\rm CI}_{\alpha}^{0}\left(Z,k,q\right),\Theta_0\left(k\right)\right)\lesssim k^{\frac{2}{q}}\frac{\log p}{n}\sigma_0^2.
\label{eq: length of CI known}
\end{equation}
In particular, for the case $q=2$, \eqref{eq: coverage of CI known} and \eqref{eq: length of CI known} also hold for the estimator $\betal$ defined in \eqref{eq: Lasso} with $A>\sqrt{2}$.
\end{Proposition}

This result shows that the confidence interval ${\rm CI}_{\alpha}^{0} \left(Z,k,q\right)$ achieves the minimax rate given in  \eqref{eq: lq minimax}.
In contrast to the $\ell_2$ loss where the two-sided confidence interval \eqref{eq: propose CI for square norm loss} is significantly shorter than the one-sided interval and achieves the optimal rate over the regime $\frac{\sqrt{n}}{\log p}\lesssim k\lesssim \frac{n}{\log p}$, for the $\ell_q$ loss with $1\leq q<2$, the one-sided confidence interval achieves the optimal rate given in \eqref{eq: lq minimax}.

We now consider adaptivity of confidence intervals. The following theorem establishes the lower bound for $\RR_{\alpha}^{*}\left(\Theta_0(k_1),\Theta_0(k_2), \widehat{\beta},\ell_{q}\right)$ with $1\leq q<2$.
\begin{Theorem}
\label{thm: adaptivity for q<2 both known}   
Suppose $0<\alpha<\frac{1}{4}$, $1\leq q<2$ and the sparsity levels $k_1,k_2$ and $k_0$ satisfy Assumption $({\rm B}2)$. 
For any  estimator $\widehat{\beta}$ satisfying  Assumption $({\rm A}2)$ with $\|\betanull\|_0\leq k_0,$ then there is some constant $c>0$ such that
\begin{equation}
\RR_{\alpha}^{*}\left(\Theta_0(k_1),\Theta_0(k_2), \widehat{\beta},\ell_{q}\right) \geq
        \begin{cases}
        c  k_2^{\frac{2}{q}} \frac{\log p}{n} \sigma_0^2 & \text{if} \quad \regimea;\\[3pt]
        c  k_2^{\frac{2}{q}-1} \frac{1}{\sqrt{n}} \sigma_0^2 & \text{if}\quad \regimeb;\\[3pt]
        c k_2^{\frac{2}{q}-{1}}k_{1} {\frac{\log p}{n}}\sigma_0^2 & \text{if}\quad \regimec.
        \end{cases}
\label{eq: lq lower bound adaptation}
\end{equation}
In particular, if $p \geq n$  and  $\betal$ is the Lasso estimator defined in \eqref{eq: Lasso} with $A> 4\sqrt{2}$, the above lower bounds can be achieved.
\end{Theorem}
The lower bounds of Theorem \ref{thm: adaptivity for q<2 both known} imply that of Theorem \ref{thm: minimax of CI for the lq loss known variance} and
both lower bounds hold for a general class of estimators satisfying Assumption $({\rm A}2)$.
However, the lower bound \eqref{eq: lq lower bound adaptation} in Theorem \ref{thm: adaptivity for q<2 both known} has a significantly different meaning from \eqref{eq: lq lower bound} in Theorem \ref{thm: minimax of CI for the lq loss known variance} where \eqref{eq: lq lower bound adaptation} quantifies the cost of adaptation without knowing the sparsity level. For the Lasso estimator $\betal$ defined in \eqref{eq: Lasso}, by comparing Theorem \ref{thm: minimax of CI for the lq loss known variance} and Theorem \ref{thm: adaptivity for q<2 both known}, we obtain $\RR_{\alpha}^{*}\left(\Theta_0(k_1),\Theta_0(k_2), \betal,\ell_{q}\right)\gg \RR_{\alpha}^{*}\left(\Theta_0(k_1), \betal,\ell_{q}\right)$ if $k_1 \ll k_2,$
which implies the impossibility of constructing adaptive confidence intervals for the case $1\leq q <2$. There exists marked difference between the case $1\leq q<2$ and the case $q=2$, where it is possible to construct adaptive confidence intervals  over the regime $\frac{\sqrt{n}}{\log p} \lesssim k \lesssim \frac{n}{\log p}$.

For the Lasso estimator $\betal$ defined in \eqref{eq: Lasso}, it is shown in Proposition \ref{prop: construction of CI known} that the confidence interval ${\rm CI}_{\alpha}^{0} \left(Z,k_2,q\right)$  defined in \eqref{eq: theoretical CI for known case} achieves the lower bound $k_2^{\frac{2}{q}} \frac{\log p}{n}\sigma_0^2$ of \eqref{eq: lq lower bound adaptation}. 
The lower bounds $k_2^{\frac{2}{q}-{1}}k_{1} {\frac{\log p}{n}}\sigma_0^2$ and $k_2^{\frac{2}{q}-1} \frac{1}{\sqrt{n}}\sigma_0^2$ of \eqref{eq: lq lower bound adaptation}  can be achieved by the following proposed confidence interval, 
\begin{equation}
{\rm CI}_{\alpha}^{2}\left(Z,k_2,q\right) = 
  \left(\left(\frac{\psi\left(Z\right)}{\frac{1}{\ns}\chi^2_{1-\frac{\alpha}{2}}\left(\ns\right)}-\sigma_0^{2}\right)_{+},\left(16 k_2\right)^{\frac{2}{q}-1}\left(\frac{\psi\left(Z\right)}{\frac{1}{\ns}\chi^2_{\frac{\alpha}{2}}\left(\ns\right)}-\sigma_0^{2}\right)_{+}\right),  
\label{eq: propose CI for q<2 norm loss}
\end{equation}
where $\psi\left(Z\right)$ is given in \eqref{eq: def of radius}. The following result verifies the above claim.
\begin{Proposition}
\label{prop: construction of CI2}
Suppose $p\geq n$, $k_1\leq k_2\lesssim \frac{n}{\log p}$ and $\betal$ is defined in \eqref{eq: Lasso} with $A>4\sqrt{2}$. 
Then ${\rm CI}_{\alpha}^{2}\left(Z,k_2,q\right)$ defined in \eqref{eq: propose CI for q<2 norm loss} satisfies, 
\begin{equation}
\liminfnp\inf_{\theta\in \Theta_0\left(k_2\right)}\PP_{\theta}\left(\|\widehat{\beta}-\beta\|_q^2\in {\rm CI}_{\alpha}^{2}\left(Z,k_2,q\right)\right)\geq 1-\alpha,
\label{eq: coverage of CI2 q<2}
\end{equation}
and
\begin{equation}
\RR\left({\rm CI}_{\alpha}^{2}\left(Z,k_2,q\right),\Theta_0\left(k_1\right)\right)\lesssim k_2^{\frac{2}{q}-1} \left(k_1\frac{\log p}{n}+\frac{1}{\sqrt{n}}\right)\sigma_0^2.
\label{eq: length of CI2 q<2}
\end{equation}
\end{Proposition}

\section{Minimaxity and adaptivity of confidence intervals over $\Theta(k)$} 
\label{sec: CI for unknown}

In this section, we focus on the case of unknown $\Sigma$ and $\sigma$ and establish the rates of convergence for the minimax expected length of confidence intervals for $\glosss$ with $1\leq q \leq 2$ over $\Theta(k)$ defined in \eqref{eq: parameter space with unknown variance}.  We also study the possibility of adaptivity of confidence intervals for $\glosss$. The following theorem establishes the lower bounds for the benchmark quantities  $\RR_{\alpha}^{*}\left(\Theta\left(k_i\right),\widehat{\beta},\ell_{q}\right)$ with $i=1,2$ and $\RR_{\alpha}^{*}\left(\Theta\left(k_1\right), \Theta\left(k_2\right),\widehat{\beta},\ell_{q}\right)$. 
\begin{Theorem}
\label{thm: non-adaptivity of CI for the global loss}
Suppose that $0<\alpha<\frac{1}{4}$, $1\leq q \leq 2$ and the sparsity levels $k_1,k_2$ and $k_0$ satisfy Assumption $({\rm B}2)$.
For any estimator $\widehat{\beta}$ satisfying Assumption $({\rm A}1)$ at $\theta_0=(\betanull,{\rm I},\sigma_0)$ with $\|\betanull\|_0\le k_0$,  there is a constant $c>0$ such that
\begin{equation}
\RR_{\alpha}^{*}\left(\Theta\left(k_i\right), \widehat{\beta},\ell_{q}\right) \geq c k_i^{\frac{2}{q}} \frac{\log p}{n}, \quad \text{for} \quad i=1,2;
\label{eq: non-adaptivity of CI for the global loss}
\end{equation}
\begin{equation}
\RR_{\alpha}^{*}\left(\left\{\theta_0\right\},\Theta\left(k_2\right), \widehat{\beta},\ell_{q}\right)\geq c k_2^{\frac{2}{q}}{\frac{\log p}{n}}.
\label{eq: non-adaptivity global loss}
\end{equation}
In particular, if $\betasl$ is the scaled Lasso estimator defined in \eqref{eq: scaled Lasso} with $A>2\sqrt{2}$, then
the above lower bounds can be achieved.
\end{Theorem}

\begin{figure}[h!]
\centering
\includegraphics[width=3in,height=1.2in]{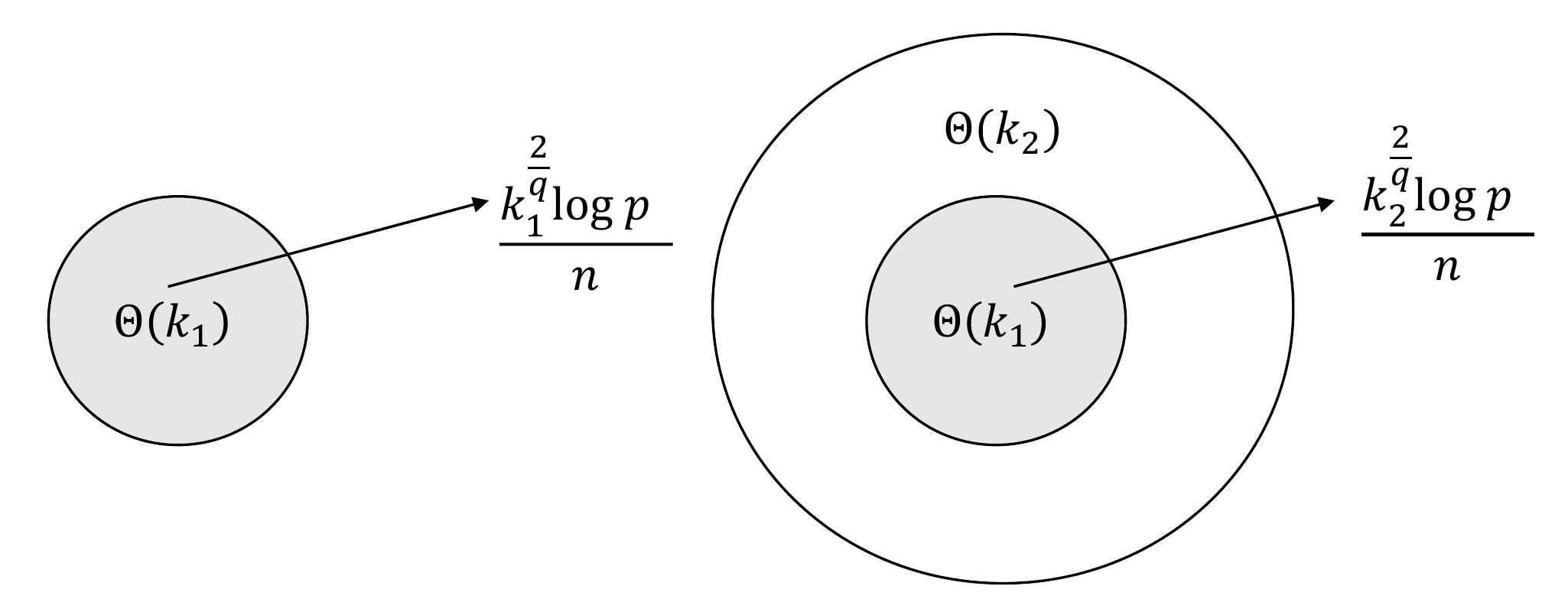}
\caption{Illustration of $\RR_{\alpha}^{*}\left(\Theta\left(k_1\right), \betasl,\ell_{q}\right)$ (left) and $\RR_{\alpha}^{*}\left(\Theta\left(k_1\right),\Theta\left(k_2\right), \betasl,\ell_{q}\right)$ (right).}
\label{fig:adaptivity behavior of lq norm unknown}
\end{figure}

The lower bounds \eqref{eq: non-adaptivity of CI for the global loss} and \eqref{eq: non-adaptivity global loss} hold for any $\widehat{\beta}$ satisfying Assumption $({\rm A}1)$ at an interior point $\theta_0$, including the scaled Lasso estimator as a special case.
We demonstrate the impossibility of adaptivity of confidence intervals for the $\ell_{q}$ loss of  the scaled Lasso estimator $\betasl$ defined in \eqref{eq: scaled Lasso}. 
Since $\RR_{\alpha}^{*}\left(\Theta\left(k_1\right),\Theta\left(k_2\right), \betasl,\ell_{q}\right)\geq \RR_{\alpha}^{*}\left(\left\{\theta_0\right\},\Theta\left(k_2\right), \betasl,\ell_{q}\right),$ by \eqref{eq: non-adaptivity global loss}, we have $\RR_{\alpha}^{*}\left(\Theta\left(k_1\right),\Theta\left(k_2\right), \betasl,\ell_{q}\right)\gg \RR_{\alpha}^{*}\left(\Theta\left(k_1\right), \betasl,\ell_{q}\right)$ if $k_1 \ll k_2.$ The comparison of $\RR_{\alpha}^{*}\left(\Theta\left(k_1\right), \betasl,\ell_{q}\right)$  and $\RR_{\alpha}^{*}\left(\Theta\left(k_1\right),\Theta\left(k_2\right), \betasl,\ell_{q}\right)$ is illustrated in Figure \ref{fig:adaptivity behavior of lq norm unknown}. Referring to the adaptivity defined in \eqref{eq: def of adaptivity for global loss}, it is impossible to construct adaptive confidence intervals for $\|\betasl-
\beta\|_q^2$. 

Theorem \ref{thm: non-adaptivity of CI for the global loss} shows that for any confidence interval $\gCIZR$ for the loss of any estimator $\widehat{\beta}$ satisfying Assumption $({\rm A}1)$, under the coverage constraint that $\gCIZR \in \II_{\alpha}\left(\Theta\left(k_2\right), \widehat{\beta},\ell_{q}\right)$, its expected length at any given $\theta_{0}=\left(\betanull ,\rm I, \sigma \right) \in  \Theta\left(k_0\right)$ must be of order
$k_2^{\frac{2}{q}}{\frac{\log p}{n}}.$
In contrast to Theorem \ref{thm: adaptation of CI for the l2 loss both known} and \ref{thm: adaptivity for q<2 both known},  Theorem \ref{thm: non-adaptivity of CI for the global loss}  demonstrates that confidence intervals must be long at a large subset of points in the parameter space, not just at a small number of ``unlucky" points. Therefore, the lack of adaptivity for confidence intervals is not due to the conservativeness of the minimax framework.

In the following, we detail the construction of confidence intervals for $\|\betasl-\beta\|_q^{2}$. 
The construction of confidence intervals is based on the following definition of restricted eigenvalue, which is introduced in \citep{bickel2009simultaneous},
\begin{equation}
\kappa(X,k,s,\alpha_0)=\min_{\substack{J_{0}\subset \{1,\cdots,p\},\\|J_0|\leq k}} \min_{\substack{\delta\neq 0,\\ \|\delta_{J_0^c}\|_1\leq \alpha_0 \|\delta_{J_0}\|_1}} \frac{\|X\delta\|_2}{\sqrt{n} \|\delta_{J_{01}}\|_2},
\end{equation}
where $J_1$ denotes the subset corresponding to the $s$ largest in absolute value coordinates of $\delta$ outside of $J_0$ and $J_{01}=J_{0} \cup J_1$.
Define the event 
$\mathcal{B}=\left\{\hat{\sigma}\leq \log p\right\}.$
The confidence interval for $\|\betasl-\beta\|_q^2$ is defined as
\begin{equation}
{\rm CI}_{\alpha}\left(Z,k,q\right)= \left\{
  \begin{array}{cl}
    \left[0, \; \varphi\left(Z,k,q\right) \right]    &\quad \text{on}\; \mathcal{B}\\
   \left\{0\right\}  &\quad \text{on}\; \mathcal{B}^{c},
  \end{array}
\right.
\label{eq: constructed CI for multi sparse}
\end{equation}
where
\begin{equation*}
\varphi\left(Z,k,q\right)=\min\left\{\left(\frac{16 A{\max \|\Xj\|_2^2} \widehat{\sigma}}{{n}\kappa^2\left(X,k,k,3\left(\frac{\max \|\Xj\|_2}{\min \|\Xj\|_2}\right) \right)}\right)^2 k^{\frac{2}{q}}\frac{\log p}{n}, \; \left(k^{\frac{2}{q}}\frac{\log p}{n} \log p\right)\widehat{\sigma}^2\right\}. 
\label{eq: scaled Lasso bound for beta}
\end{equation*}

Properties 
of ${\rm CI}_{\alpha}\left(Z,k,q\right)$ are established as follows.
\begin{Proposition}
\label{prop: construction of CI unknown}
Suppose $k \lesssim \frac{n}{\log p}$ and $\betasl$ is the estimator defined in \eqref{eq: scaled Lasso} with $A>2\sqrt{2}$. 
For $1\leq q \leq 2$, then ${\rm CI}_{\alpha}\left(Z,k,q\right)$ defined in \eqref{eq: constructed CI for multi sparse} satisfies the following properties,
\begin{equation}
\liminfnp \inf_{\theta\in \Theta\left(k\right)}\PP_{\theta}\left(\|\widehat{\beta}-\beta\|_q^2\in {\rm CI}_{\alpha}\left(Z,k,q\right)\right) =1, 
\label{eq: coverage of CI unknown}
\end{equation}
and
\begin{equation}
\RR\left({\rm CI}_{\alpha}\left(Z,k,q\right),\Theta\left(k\right)\right)\lesssim k^{\frac{2}{q}}\frac{\log p}{n}.
\label{eq: length of CI unknown}
\end{equation}
\end{Proposition}
Proposition \ref{prop: construction of CI unknown} shows that the confidence interval ${\rm CI}_{\alpha}\left(Z,k_i,q\right)$ defined in \eqref{eq: constructed CI for multi sparse} achieves the lower bound
in \eqref{eq: non-adaptivity of CI for the global loss}, for $i=1,2$,
and the confidence interval ${\rm CI}_{\alpha}\left(Z,k_2,q\right)$ defined in \eqref{eq: constructed CI for multi sparse} achieves the lower bound in \eqref{eq: non-adaptivity global loss}.

\section{Estimation of the $\ell_q$ loss of rate-optimal estimators}
\label{sec: rate optimal estimators}

We have established the minimax lower bounds for the estimation accuracy of the loss of a broad class of estimators $\widehat{\beta}$ satisfying the weak assumptions (A1) or (A2) and also demonstrated that such minimax lower bounds are sharp for the Lasso estimator or scaled Lasso estimator. In this section, we will show that the minimax lower bounds are sharp for the class of rate-optimal estimators satisfying the following Assumption (A).

\begin{enumerate} 
\item[(A)] The estimator $\widehat{\beta}$ satisfies, for $k\ll \frac{n}{\log p}$, 
\begin{equation}
\sup_{\theta\in \Theta(k)}\PP_{\theta} \left(\|\widehat{\beta}-\beta\|_{q}^2\geq C^* \|\beta\|_0^{\frac{2}{q}} {\frac{\log p}{n}}\right)\leq C p^{-\delta},
 \label{eq: strong adaptivity assumption}
 \end{equation}
for constants $\delta>0$, $C^*>0$ and $C>0$.
 \end{enumerate}
 We say an estimator $\widehat{\beta}$ is  rate-optimal if it satisfies Assumption (A).
 As shown in \citep{candes2007dantzig,bickel2009simultaneous,belloni2011square,sun2012scaled}, Lasso, Dantzig Selector, scaled Lasso and square-root Lasso are rate-optimal when the tuning parameter is chosen properly.
We shall stress that Assumption (A) implies  Assumptions (A1) and (A2). Assumption (A) requires the estimator $\widehat{\beta}$ to perform well over the whole parameter space $\Theta(k)$ while Assumptions (A1) and (A2) only require $\widehat{\beta}$ to perform well at a single point or over a proper subset. The following proposition shows that the minimax lower bounds established in Theorem \ref{thm: estimating the global loss known} to Theorem \ref{thm: non-adaptivity of CI for the global loss} can be achieved for the class of rate-optimal estimators.


\begin{Proposition}
\label{prop: matched upper bound}
Let $\widehat{\beta}$ be an estimator satisfying Assumption {\rm (A)}.
\begin{enumerate}
\item There exist estimators of the loss $\|\widehat{\beta}-\beta\|_{q}^2$ with $1\leq q<2$ achieving, up to a constant factor, the minimax lower bounds \eqref{eq: key lower bound in high prob known application general q} in Theorem \ref{thm: estimating the global loss known} and  \eqref{eq: lq lower bound} in Theorem \ref{thm: minimax of CI for the lq loss known variance}  and estimators of loss  $\|\widehat{\beta}-\beta\|_{q}^2$ with $1\leq q\leq 2$ achieving, up to a constant factor, the minimax lower bounds  \eqref{eq: key lower bound in high prob unknown application} in Theorem \ref{thm: estimating the global loss} and \eqref{eq: non-adaptivity of CI for the global loss} in Theorem \ref{thm: non-adaptivity of CI for the global loss}. 

\item Suppose that the estimator $\widehat{\beta}$ is constructed based on the subsample $Z^{\left(1\right)}=\left(y^{\left(1\right)},X^{\left(1\right)}\right)$, then there exist estimators of the loss $\|\widehat{\beta}-\beta\|_{2}^2$ achieving, up to a constant factor, the minimax lower bounds \eqref{eq: key lower bound in high prob known application q=2} in Theorem \ref{thm: estimating the global loss known}, \eqref{eq: l2 lower bound} in Theorem \ref{thm: minimax of CI for the l2 loss both known} and \eqref{eq: l2 lower bound adaptation} in Theorem \ref{thm: adaptation of CI for the l2 loss both known}.

\item Suppose the estimator $\widehat{\beta}$ is constructed based on the subsample $Z^{\left(1\right)}=\left(y^{\left(1\right)},X^{\left(1\right)}\right)$ and it satisfies Assumption (A) with $\delta>2$ and the assumption $\|(\widehat{\beta}-\beta)_{S^{c}}\|_1\leq c^{*} \|(\widehat{\beta}-\beta)_{S}\|_1$ where $S={\rm supp}(\beta)$. Then for $p \ge n$ there exist estimators of the loss $\|\widehat{\beta}-\beta\|_{q}^2$ with $1\leq q<2$ achieving the minimax lower bounds \eqref{eq: lq lower bound adaptation} in Theorem \ref{thm: adaptivity for q<2 both known}. 
\end{enumerate}
\end{Proposition}
For reasons of space, we do not discuss the detailed construction of confidence intervals achieving these minimax lower bounds here and postpone the construction to the proof of Proposition \ref{prop: matched upper bound}. 
\begin{Remark} \rm
Sample splitting has been widely used in the literature. For example, the condition that $\widehat{\beta}$ is constructed based on the subsample $Z^{\left(1\right)}=\left(y^{\left(1\right)},X^{\left(1\right)}\right)$ has been introduced in \citep{nickl2013confidence} for constructing confidence sets for $\beta$ and in \citep{janson2015eigenprism} for constructing confidence intervals for the $\ell_2$ loss. Such a condition is imposed purely  for technical reasons to create the independence between the estimator $\widehat{\beta}$ and the subsample $Z^{\left(2\right)}=\left(y^{\left(2\right)},X^{\left(2\right)}\right)$, which is used to evaluate the $\ell_{q}$ loss of the estimator $\widehat{\beta}$. 
As shown in \citep{bickel2009simultaneous}, the assumption $\|(\widehat{\beta}-\beta)_{S^{c}}\|_1\leq c^{*} \|(\widehat{\beta}-\beta)_{S}\|_1$ is satisfied for Lasso and Dantzig Selector.
\end{Remark}

\section{General tools for minimax lower bounds}
\label{sec: lower bound tool}

A major step in our analysis is to establish rate sharp lower bounds for the estimation error and  expected length of confidence intervals for the $\ell_q$ loss. We introduce in this section new technical tools that are needed to establish these lower bounds. 

A significant distinction of the lower bound results given in the previous sections from those for the traditional parameter estimation problems is that the constraint is on the performance of the estimator $\widehat{\beta}$ of the regression vector $\beta$, but the lower bounds are on the difficulty of estimating its loss $\|\widehat{\beta} - \beta\|_q^2$.  It is necessary to develop new lower bound techniques 
 to establish rate-optimal lower bounds for the estimation error and the expected length of confidence intervals for the loss $\|\widehat{\beta} - \beta\|_q^2$. These technical tools may also be of independent interest. 

We begin with notation. Let $Z$ denote a random variable whose distribution is indexed by some parameter $\theta \in \Theta$ and let $\pi$ denote a prior on the parameter space $\Theta$. We will use $f_{\theta}(z)$ to denote the density of $Z$ given $\theta$ and $f_{\pi}\left(z\right)$ to denote the marginal density of $Z$ under the prior $\pi$. Let  $\PP_{\pi}$ denote the distribution of $Z$ corresponding to $f_{\pi}\left(z\right)$, i.e., $\PP_{\pi}\left(\mathcal{A}\right)=\int 1_{z\in \mathcal{A}} f_{\pi}\left(z\right) dz,$ where $1_{z\in \mathcal{A}}$ is the indicator function.
For a function $g$, we write $\E_{\pi}\left(g(Z)\right)$ for the expectation under $f_{\pi}$. More specifically, 
$f_{\pi}\left(z\right)=\int f_{\theta}\left(z\right) \pi\left(\theta\right) d\theta$ and $\E_{\pi}\left(g(Z)\right)=\int g\left(z\right) f_{\pi}\left(z\right) dz.$
The $L_1$ distance between two probability distributions with densities $f_0$ and $f_1$ is given by
$\TV(f_{1},f_{0})=\int \left|f_1(z)-f_0(z)\right| dz.$ The following theorem establishes the minimax lower bounds for the estimation error and the expected length of confidence intervals for the $\ell_{q}$ loss, under the constraint that $\widehat{\beta}$ is a good estimator at at least one interior point.
\begin{Theorem}
\label{thm: lower bound loss estimation generalization}
Suppose $0<\alpha, \alpha_0<\frac{1}{4}$, $1\leq q\leq 2$, $\Sigma_{0}$ is positive definite, $\theta_0=\left(\betanull, { \Sigma_0},\sigma_0 \right) \in \Theta$, and $\FF \subset \Theta$.
 Define $\distalter=\min_{\theta \in \FF} \|\beta\left(\theta\right)-\betanull\|_{q}$. Let $\pi$ denote a prior over the parameter space $\FF$.  If an estimator $\widehat{\beta}$ satisfies 
\begin{equation}
\PP_{\theta_0}\left(\|\widehat{\beta}-\betanull\|_{q}^2 \leq \frac{1}{16} \distalter^2\right)\geq 1-\alpha_0,
\label{eq: adaptive risk assumption high prob}
\end{equation}
then
\begin{equation}
\inf_{\widehat{L}_{q}} \sup_{\theta \in \left\{\theta_0\right\}\cup\FF} \PP_{\theta}\left(|\widehat{L}_{q}-\glosss| \geq \frac{1}{4} d^2\right) \geq \bar{c}_1,
\label{eq: key lower bound in high prob unknown}
\end{equation}
and
\begin{equation}
\RR_{\alpha}^{*}\left(\left\{\theta_0\right\},\Theta, \widehat{\beta},\ell_{q}\right)=\inf_{{\rm CI}_{\alpha}\left(\widehat{\beta},\ell_{q},Z\right)\in \II_{\alpha}\left(\Theta,\widehat{\beta}, \ell_{q}\right)} \E_{\theta_0} \RR\left({\rm CI}_{\alpha}\left(\widehat{\beta},\ell_{q}, Z\right)\right)\geq c_2^* d^2,
\label{eq: lower bound for CI of loss}\end{equation}
where
 $\bar{c}_1=\min\left\{\frac{1}{10}, \left(\frac{9}{10}-\alpha_0-\TV\left(f_{\pi},f_{\theta_0}\right) \right)_{+}\right\}$ and  $c_{2}^{*}=\frac{1}{2}\left(1-2\alpha-\alpha_0-2\TV\left(f_{\pi},f_{\theta_0}\right)\right)_{+}.$
\end{Theorem}

\begin{Remark} \rm
The minimax lower bound \eqref{eq: key lower bound in high prob unknown} for the estimation error and \eqref{eq: lower bound for CI of loss} for the expected length of confidence intervals hold as long as the estimator $\widehat{\beta}$ estimates $\beta$ well at an interior point $\theta_0$. Besides Condition \eqref{eq: adaptive risk assumption high prob}, another key ingredient for the lower bounds  \eqref{eq: key lower bound in high prob unknown} and \eqref{eq: lower bound for CI of loss} is to construct the least favorable space $\FF$ with the prior $\pi$ such that the marginal distributions $f_{\pi}$ and $f_{\theta_0}$ are non-distinguishable. 
For the estimation lower bound \eqref{eq: key lower bound in high prob unknown}, constraining that $\|\widehat{\beta}-\betanull\|_{q}^2$  can be well estimated at $\theta_0$, due to the non-distinguishability between $f_{\pi}$ and $f_{\theta_0}$, we can establish that the loss $\glosss$ cannot be estimated well over $\FF$. For the lower bound \eqref{eq: lower bound for CI of loss}, by Condition \eqref{eq: adaptive risk assumption high prob} and the non-distinguishability between $f_{\pi}$ and $f_{\theta_0}$, we will show that $\glosss$ over $\FF$ is much larger than $\|\widehat{\beta}-\betanull\|_{q}^2$ and hence the honest confidence intervals must be sufficiently long.
\end{Remark}

Theorem \ref{thm: lower bound loss estimation generalization} is used to establish the minimax lower bounds for both the estimation error and the expected length of confidence intervals of the $\ell_{q}$ loss over $\Theta(k)$. By taking $\theta_0 \in \Theta(k_0)$ and $\Theta=\Theta(k)$, Theorem \ref{thm: estimating the global loss} follows from \eqref{eq: key lower bound in high prob unknown} with a properly constructed subset $\FF\subset \Theta(k)$.  
By taking $\theta_0\in \Theta(k_0)$ and $\Theta=\Theta(k_2)$, Theorem \ref{thm: non-adaptivity of CI for the global loss} follows from \eqref{eq: lower bound for CI of loss}  with a properly constructed $\FF\subset\Theta(k_2)$.  In both cases, Assumption $({\rm A}1)$ implies Condition \eqref{eq: adaptive risk assumption high prob}.

Several minimax lower bounds over $\Theta_0(k)$ can also be implied by Theorem \ref{thm: lower bound loss estimation generalization}. For the estimation error, the minimax lower bounds \eqref{eq: key lower bound in high prob known application q=2}  and \eqref{eq: key lower bound in high prob known application general q} over the regime $k\lesssim \frac{\sqrt{n}}{\log p}$ in Theorem \ref{thm: estimating the global loss known} follow from \eqref{eq: key lower bound in high prob unknown}. For the expected length of confidence intervals, the minimax lower bounds \eqref{eq: l2 lower bound adaptation} in Theorem \ref{thm: adaptation of CI for the l2 loss both known}  and \eqref{eq: lq lower bound adaptation} in the regions $\regimea$ and $\regimeb$ in Theorem \ref{thm: adaptivity for q<2 both known} follow from \eqref{eq: lower bound for CI of loss}. In these cases, Assumption $({\rm A}1)$ or $({\rm A}2)$ can guarantee that Condition \eqref{eq: adaptive risk assumption high prob} is satisfied.
However, the minimax lower bound for estimation error \eqref{eq: key lower bound in high prob known application general q} in the region $\frac{\sqrt{n}}{\log p}\leq k\lesssim \frac{{n}}{\log p}$ and for the expected length of confidence intervals \eqref{eq: lq lower bound adaptation} in the region $\regimec$  cannot be established using the above theorem. The following theorem, which requires testing a composite null against a composite alternative,  establishes the refined minimax lower bounds over $\Theta_0(k)$.  

\begin{Theorem}
\label{thm: lower bound for CI of risk high prob known variance}
Let $0<\alpha,\alpha_0<\frac{1}{4}$, $1\leq q\leq 2$, and $\theta_0=\left(\betanull,{\Sigma_0},\sigma_0\right)$ where {$\Sigma_{0}$ is a positive definite matrix}. 
Let $k_1$ and $k_2$ be two sparsity levels.
Assume that for $i=1,2$ there exist parameter spaces $\FF_i\subset \left\{\left(\beta,\Sigma_0,\sigma_0 \right): \|\beta\|_0\leq k_i\right\}$  such that for given $\disttwonorm_i$ and $\distalter_i$  
\[
\left\|{\Sigma_0}\left(\beta\left(\theta\right)-\betanull\right) \right\|_2=\disttwonorm_i \quad \text{\rm and} \quad \|\beta\left(\theta\right)-\betanull\|_{q}=\distalter_i, \quad\mbox{\rm for all $\theta \in \FF_i$.}
\]
Let $\pi_i$ denote a prior over the parameter space $\FF_i$ for $i=1,2$. Suppose that for $\theta_1=\left(\betanull,{\Sigma_0}, \sigma_0^2+\disttwonorm_1^2 \right)$ and $\theta_2=\left(\betanull,{\Sigma_0}, \sigma_0^2+\disttwonorm_2^2\right)$, there exist constants $c_1, c_2 >0$ such that 
\begin{equation}
\PP_{\theta_i}\left(\|\widehat{\beta}-\betanull\|_{q}^2\leq c^2_i \distalter_i^2\right)\geq 1-\alpha_0, \quad \text{for} \quad i=1,2.
\label{eq: adaptive high prob risk at zero general}
\end{equation}
Then we have
\begin{equation}
\inf_{\widehat{L}_{q}} \sup_{\theta \in \FF_1\cup\FF_2} \PP_{\theta}\left(|\widehat{L}_{q}-\glosss| \geq c^{*}_3 \distalter_2^2\right) \geq \bar{c}_3,
\label{eq: key lower bound in high prob known}
\end{equation}
and
\begin{equation}
\RR_{\alpha}^{*}\left(\Theta_0\left(k_1\right),\Theta_0\left(k_2\right),\widehat{\beta},\ell_{q}\right)\\
\geq  c_{4}^{*} \left(\left(1-c_2\right)^2 \distalter_2^2-\left(1+c_1\right)^2\distalter_1^2\right)_{+},
\label{eq: general non-adaptive lower bound}
\end{equation}
where
{\footnotesize $c^{*}_3=\min\left\{\frac{1}{4} ,\left(\left(1-c_2\right)^2 -\frac{1}{4}-(1+c_1)^2\frac{\distalter_1^2}{\distalter_2^2}\right)_{+}\right\},$ $c_{4}^{*}=\left(1-2\alpha_0-2\alpha-\sum_{i=1}^{2}\TV\left(f_{\pi_i},f_{\theta_i}\right)-2\TV\left(f_{\pi_2},f_{\pi_1}\right)\right)_{+}$ and $\bar{c}_3=\min\left\{\frac{1}{10}, \left(\frac{9}{10}-2\alpha_0-\sum_{i=1}^{2}\TV\left(f_{\pi_i},f_{\theta_i}\right)-2\TV\left(f_{\pi_2},f_{\pi_1}\right) \right)_{+}\right\}.$}
\end{Theorem}

\begin{Remark}\rm
As long as the estimator $\widehat{\beta}$ performs well at two points, $\theta_1$ and $\theta_2$, the minimax lower bounds \eqref{eq: key lower bound in high prob known} for the estimation error and \eqref{eq: general non-adaptive lower bound} for the expected length of confidence intervals hold. Note that $\theta_i$ in the above theorem does not belong to the parameter space $ \left\{\left(\beta,\Sigma_0,\sigma_0 \right): \|\beta\|_0\leq k_i\right\}$, for $i=1,2$. 
In contrast to Theorem \ref{thm: lower bound loss estimation generalization}, Theorem \ref{thm: lower bound for CI of risk high prob known variance} compares composite hypotheses $\FF_1$ and $\FF_2$, which will lead to a sharper lower bound than comparing the simple null $\{\theta_0\}$ with the composite alternative $\FF$. 
For simplicity, we construct least favorable parameter spaces $\FF_i$ such that the points in $\FF_i$ is of fixed $\ell_2$ distance and fixed $\ell_{q}$ distance to $\betanull$, for $i=1,2$, respectively. More importantly, we construct $\FF_1$ with the prior $\pi_1$ and $\FF_2$ with the prior $\pi_2$ such that $f_{\pi_1}$ and $f_{\pi_2}$ are not distinguishable, where $\theta_1$ and $\theta_2$ are introduced to facilitate the comparison.  By Condition \eqref{eq: adaptive high prob risk at zero general} and the construction of $\FF_1$ and $\FF_2$, we establish that the $\ell_{q}$ loss cannot be simultaneously estimated well over $\FF_1$ and $\FF_2$. For the lower bound \eqref{eq: general non-adaptive lower bound}, under the same conditions, it is shown that the $\ell_{q}$ loss over $\FF_1$ and $\FF_2$ are far apart and any confidence interval with guaranteed coverage probability over $\FF_1\cup \FF_2$ must be sufficiently long. Due to the prior information  $\Sigma={\rm I}$ and $\sigma=\sigma_0$, the lower bound construction over $\Theta_0(k)$ is more involved than that over $\Theta(k)$. We shall stress that the construction of $\FF_1$ and $\FF_2$ and the comparison between composite hypotheses are of independent interest.
\end{Remark}
The minimax lower bound \eqref{eq: key lower bound in high prob known application general q} in the region $\frac{\sqrt{n}}{\log p}\lesssim k\lesssim \frac{{n}}{\log p}$ follows from \eqref{eq: key lower bound in high prob known} and the minimax lower bound \eqref{eq: lq lower bound adaptation} in the region $\regimec$ for the expected length of confidence intervals follows from \eqref{eq: general non-adaptive lower bound}. In these cases, $\Sigma_0$ is taken as $\rm I$ and Assumption $({\rm A}2)$ implies Condition \eqref{eq: adaptive high prob risk at zero general}.
\section{An intermediate setting with known $\sigma=\sigma_0$ and unknown $\Sigma$}
\label{sec: more general parameter space}

The results given in Sections \ref{sec: CI for known} and \ref{sec: CI for unknown} show the significant difference between $\Theta_0(k)$ and $\Theta(k)$ in terms of minimaxity and adaptivity of confidence intervals for $\|\widehat{\beta}-\beta\|_q^2$. 
$\Theta_0(k)$ is for the simple setting with known design covariance matrix $\Sigma={\rm I}$ and known noise level $\sigma=\sigma_0$, and the $\Theta(k)$ is for unknown $\Sigma$ and $\sigma$.
In this section, we further consider minimaxity and adaptivity of confidence intervals for $\|\widehat{\beta}-\beta\|_q^2$  in an intermediate setting where  the noise level $\sigma=\sigma_0$ is known and $\Sigma$ is unknown but of certain structure. Specifically, we consider the following parameter space, 
\begin{equation}
\Theta_{\sigma_0}(k,s)=\left\{\left(\beta,\Sigma,\sigma_0\right): \begin{aligned}
&\|\beta\|_0\leq k, \; \;\frac{1}{M_1} \leq \lambda_{\min}\left(\Sigma\right)\leq \lambda_{\max}\left(\Sigma\right) \leq M_1\\
&\|\Sigma^{-1}\|_{L_1}\leq M, \quad \max_{1\leq i\leq p}\|\left(\Sigma^{-1}\right)_{i\cdot}\|_0\leq s
\end{aligned}
\right\},
\end{equation}
for some constants $M_1 \ge 1$ and $M>0$. $\Theta_{\sigma_0}(k,s)$ basically assumes known noise level $\sigma$ and imposes sparsity conditions on the precision matrix of the random design. This parameter space is similar to those used in the literature of sparse linear regression with random design \citep{van2014asymptotically,chernozhukov2015post,chernozhukov2015valid}. $\Theta_{\sigma_0}(k,s)$ has two sparsity parameters where $k$ represents the sparsity of $\beta$ and $s$ represents the maximum row sparsity of the precision matrix $\Sigma^{-1}$. Note that $\Theta_0(k)\subset \Theta_{\sigma_0}(k,s) \subset \Theta(k)$ and $\Theta_0(k)$  is a special case of $\Theta_{\sigma_0}(k,s)$ with $M_1=1$. 

Under the assumption $s\ll \sqrt{n/\log p}$, the minimaxity and adaptivity lower bounds for the expected length of confidence intervals for $\|\widehat{\beta}-\beta\|_q^2$ with $1\leq q<2$ over $\Theta_{\sigma_0}(k,s)$ are the same as those over $\Theta_0(k)$. That is, Theorems \ref{thm: minimax of CI for the lq loss known variance} and \ref{thm: adaptivity for q<2 both known} hold with $\Theta_{0}(k_1)$, $\Theta_{0}(k_2)$, and $\Theta_{0}(k)$ replaced by $\Theta_{\sigma_0}(k_1,s)$, $\Theta_{\sigma_0}(k_2,s)$, and $\Theta_{\sigma_0}(k,s)$, respectively.  For the case $q=2$, the following theorem establishes the minimaxity and adaptivity lower bounds for the expected length of confidence intervals for $\|\widehat{\beta}-\beta\|_2^2$ over $\Theta_{\sigma_0}(k,s).$
\begin{Theorem}
\label{thm: adaptivity for known variance and sparse precision with q=2}
Suppose $0<\alpha,\alpha_0<{1}/{4}$, $M_1>1$, $s \ll \sqrt{{n}/{\log p}}$ and the sparsity levels $k_1,k_2$ and $k_0$ satisfy Assumption $({\rm B}2)$ with the constant $c_0$ replaced by $c_0^{*}$ defined in \eqref{eq: key constant}. 
For any  estimator $\widehat{\beta}$ satisfying  
\begin{equation}
\sup_{\theta \in \Theta(k_0)}\PP_{\theta}\left(\|\widehat{\beta}-\betanull\|_q^2\geq C^{*} \|\betanull\|_0^{\frac{2}{q}} {\frac{\log p}{n}} \sigma^2 \right)\leq \alpha_0,
\label{eq: adaptation high probability lq difference more general}
\end{equation}
with a constant $C^{*}>0$, then there is some constant $c>0$ such that
\begin{equation}
\begin{aligned}
\RR_{\alpha}^{*}\left(\Theta_{\sigma_0}(k_1,s),\Theta_{\sigma_0}(k_2,s), \widehat{\beta},\ell_{2}\right) \geq c \min\left\{k_2 \frac{\log p}{n}, \max\left\{k_1{\frac{\log p}{n}},\frac{1}{\sqrt{n}}\right\}\right\} \sigma_0^2
\end{aligned}
\label{eq: l2 lower bound adaptation general}
\end{equation}
and
\begin{equation}
\RR_{\alpha}^{*}\left(\Theta_{\sigma_0}(k_i,s), \widehat{\beta},\ell_{2}\right) \geq c \frac{k_i \log p}{n} \sigma_0^2 \quad \text{and} \quad i=1,2.
\label{eq: l2 minimax rate general}
\end{equation}
In particular, if $p \geq n$  and  $\widehat{\beta}$ is constructed based on the subsample $Z^{\left(1\right)}=\left(y^{\left(1\right)},X^{\left(1\right)}\right)$ and satisfies Assumption (A) with $\delta>2$, the above lower bounds can be attained. 
\end{Theorem}
In contrast to Theorems \ref{thm: minimax of CI for the l2 loss both known} and \ref{thm: adaptation of CI for the l2 loss both known}, the lower bounds for the case $q=2$ change in the absence of the prior knowledge $\Sigma={\rm I}$ but the possibility of adaptivity of confidence intervals over $\Theta_{\sigma_0}(k,s)$ is similar to that over $\Theta_0(k)$. 
Since the Lasso estimator $\betal$ defined in \eqref{eq: Lasso} with $A> 4\sqrt{2}$ satisfies Assumption (A) with $\delta>2$, by Theorem \ref{thm: adaptivity for known variance and sparse precision with q=2}, the minimax lower bounds \eqref{eq: l2 lower bound adaptation general} and \eqref{eq: l2 minimax rate general} can be attained for $\betal$.
For $\betal$, only when $\frac{\sqrt{n}}{\log p} \lesssim k_1\leq k_2\lesssim \frac{n}{\log p}$,  $\RR_{\alpha}^{*}\left(\Theta_{\sigma_0}(k_1,s), {\betal},\ell_{2}\right)\asymp \RR_{\alpha}^{*}\left(\Theta_{\sigma_0}(k_1,s),\Theta_{\sigma_0}(k_2,s), \widehat{\beta},\ell_{2}\right)\asymp \frac{k_1 \log p}{n}$ and adaptation between $\Theta_{\sigma_0}(k_1,s)$ and $\Theta_{\sigma_0}(k_2,s)$ is possible.
In other regimes, if $k_1\ll k_2$, then $\RR_{\alpha}^{*}\left(\Theta_{\sigma_0}(k_1,s), {\betal},\ell_{2}\right)\ll \RR_{\alpha}^{*}\left(\Theta_{\sigma_0}(k_1,s),\Theta_{\sigma_0}(k_2,s), \widehat{\beta},\ell_{2}\right)$ and adaptation between $\Theta_{\sigma_0}(k_1,s)$ and $\Theta_{\sigma_0}(k_2,s)$ is impossible. For reasons of space, more discussion on $\Theta_{\sigma_0}(k,s)$, including the construction of adaptive confidence intervals over the regime $\frac{\sqrt{n}}{\log p} \lesssim k_1\leq k_2\lesssim \frac{n}{\log p}$, is postponed to the supplement \citep{cai2016asupplement}.

\section{Minimax lower bounds for estimating $\|\beta\|_q^2$ with $1\leq q\leq 2$}
\label{sec: lq norm functional}

The lower bounds developed in this paper have broader implications. In particular, the established results imply the minimax lower bounds for  estimating $\|\beta\|_q^2$ and the expected length of confidence intervals for $\|\beta\|_q^2$ with $1\leq q\leq 2$. To build the connection, it is sufficient to note that the trivial estimator $\widehat{\beta}=0$ satisfies Assumptions (A1) and (A2) with $\betanull=0$. Then we can apply the lower bounds \eqref{eq: key lower bound in high prob known application q=2}, \eqref{eq: key lower bound in high prob known application general q} and \eqref{eq: key lower bound in high prob unknown application} to the estimator $\widehat{\beta}=0$  and establish the minimax lower bounds of estimating $\|\beta\|_q^2$,
\begin{equation}
\inf_{\widehat{L}_2} \sup_{\theta \in \Theta_0(k)} \PP_{\theta}\left(|\widehat{L}_2 - \|\beta\|_2^2| \geq c \min\left\{k{\frac{\log p}{n}},\frac{1}{\sqrt{n}}\right\}\sigma_0^2\right) \geq \delta;
\label{eq: deduced l2 norm}
\end{equation}
\begin{equation}
\inf_{\lest} \sup_{\theta \in \Theta_0(k)} \PP_{\theta}\left(|\lest- \|\beta\|_q^2| \geq c k^{\frac{2}{q}}{\frac{\log p}{n}}\sigma_0^2\right) \geq \delta, \quad \text{for} \; 1\leq q<2,
\label{eq: deduced lq norm}
\end{equation}
\begin{equation}
\inf_{\lest} \sup_{\theta \in \Theta(k)} \PP_{\theta}\left(|\lest-\|\beta\|_q^2| \geq c k^{\frac{2}{q}} {\frac{\log p}{n}}\right) \geq \delta, \quad \text{for} \; 1\leq q\leq 2,
\label{eq: deduced lq norm unknown}
\end{equation}
for some constants $\delta>0$ and $c>0$. Similarly, all the lower bounds for the expected length of confidence intervals for $\|\widehat{\beta}-\beta\|_q^2$ established in Theorem \ref{thm: minimax of CI for the l2 loss both known} to Theorem \ref{thm: non-adaptivity of CI for the global loss}  imply corresponding lower bounds for $\|\beta\|_{q}^2$.
The lower bound $\min\{k{\frac{\log p}{n}},\frac{1}{\sqrt{n}}\}\sigma_0^2$ in \eqref{eq: deduced l2 norm} is the same as the detection boundary in the sparse linear regression for the case $\Sigma={\rm I}$ and $\sigma=1$; See \citep{ingster2010detection} and \citep{arias2011global} for more details. Estimation of $\|\beta\|_2^2$ in high-dimensional linear regression has been considered in \citep{guo2016optimal} under the general setting where $\Sigma$ and $\sigma$ are unknown and the lower bound \eqref{eq: deduced lq norm unknown} with $q=2$ leads to one key component of the lower bound $c k{\frac{\log p}{n}}$ for estimating $\|\beta\|_2^2$.

\section{Proofs}
\label{sec: proof}

In this section, we present the proofs of the lower bound results. 
In Section \ref{sec: general unknown}, we establish the general lower bound result, Theorem \ref{thm: lower bound loss estimation generalization}. By applying Theorem \ref{thm: lower bound loss estimation generalization} and Theorem \ref{thm: lower bound for CI of risk high prob known variance}, we establish 
Theorems \ref{thm: adaptation of CI for the l2 loss both known} and \ref{thm: adaptivity for q<2 both known} in Section \ref{sec: known variance lower bound proof}. 
For reasons of space, the proofs of Theorems \ref{thm: estimating the global loss known}, \ref{thm: estimating the global loss}, \ref{thm: minimax of CI for the l2 loss both known}, \ref{thm: minimax of CI for the lq loss known variance}, \ref{thm: non-adaptivity of CI for the global loss}, \ref{thm: lower bound for CI of risk high prob known variance}, \ref{thm: adaptivity for known variance and sparse precision with q=2}, the upper bound results, including Propositions \ref{prop: positive estimation q=2},  \ref{prop: l2 positive CI}, \ref{prop: construction of CI known}, \ref{prop: construction of CI2}, \ref{prop: construction of CI unknown}, \ref{prop: matched upper bound} and the proofs of technical lemmas are postponed to the supplement \citep{cai2016asupplement}.

We define the $\chi^2$ distance between two density functions $f_{1}$ and $f_{0}$  by
$\chi^2(f_{1},f_{0})=\int \frac{\left(f_1(z)-f_0(z)\right)^2}{f_0(z)} dz=\int \frac{f^2_{1}(z)}{f_{0}(z)} dz-1,$
and 
it is well known that 
\begin{equation}
\TV(f_{1},f_{0})\leq \sqrt{\chi^2(f_{1},f_{0})}.
\label{eq: relation between chisq and TV}
\end{equation}
Let $\PP_{Z,\theta\sim \pi}$ denote the joint probability of $Z$ and $\theta$ with the joint density function $f(\theta,z)=f_{\theta}\left(z\right) \pi\left(\theta\right).$ We introduce the following lemma, which is used in the proofs of Theorem \ref{thm: lower bound loss estimation generalization} and Theorem \ref{thm: lower bound for CI of risk high prob known variance}.  The proof of this lemma can be found in the supplement \citep{cai2016asupplement}.
\begin{Lemma}
\label{lem: joint prob and marginal prob}
For any event $\mathcal{A}$, we have 
\begin{equation}
\PP_{\pi}\left(Z\in \mathcal{A}\right)=\PP_{Z,\theta \sim \pi}\left(Z\in \mathcal{A}\right),
\label{eq: joint prob and marginal prob}
\end{equation}
\begin{equation}
\left|\PP_{\pi_1}\left(Z \in \mathcal{A}\right)-\PP_{\pi_2}\left(Z \in \mathcal{A}\right)\right|\leq \TV\left(f_{\pi_2},f_{\pi_1}\right).
\label{eq: bound by TV}
\end{equation}
\end{Lemma}
We will write $\PP_{\pi}(\mathcal{A})$ and $\PP_{Z,\theta \sim \pi}(\mathcal{A})$ for $\PP_{\pi}(Z \in \mathcal{A})$ and $\PP_{Z,\theta \sim \pi}(Z \in \mathcal{A})$ respectively.  Recall that $\widehat{L}_{q}(Z)$ denotes a data-dependent loss estimator and $\beta(\theta)$ denotes the corresponding $\beta$ of the parameter $\theta$.
\subsection{Proof of Theorem \ref{thm: lower bound loss estimation generalization}} 
\label{sec: general unknown}
We set $c_0=\frac{1}{4}$ and $\alpha_1=\frac{1}{10}$.\\
\underline{\textbf{Proof of \eqref{eq: key lower bound in high prob unknown}}}\\
We assume 
\begin{equation}
\PP_{\theta_0} \left(\left|\widehat{L}_{q}(Z)-\lossnullr\right|\leq \frac{1}{4} \distalter^2 \right)\geq 1-\alpha_1.
\label{eq: good estimates of risk at theta0}
\end{equation}
Otherwise, we have 
\begin{equation}
\PP_{\theta_0} \left(\left|\widehat{L}_{q}(Z)-\lossnullr \right|\geq \frac{1}{4} \distalter^2 \right)\geq \alpha_1.
\label{eq: bad estimate of risk at theta0 in probability}
\end{equation}
and hence \eqref{eq: key lower bound in high prob unknown} follows.
Define the event
\begin{equation}
\mathcal{A}_0=\left\{ z:
{\lossnull}\leq c^2_0 \distalter^2\; ,\; \left|\widehat{L}_{q}(z)-\lossnull\right|\leq \frac{1}{4} \distalter^2
\right\}.
\label{eq: high prob event at null hypothesis}
\end{equation}
By \eqref{eq: adaptive risk assumption high prob}  and \eqref{eq: good estimates of risk at theta0}, we have
$\PP_{\theta_0}\left(\mathcal{A}_0\right)\geq 1-\alpha_0-\alpha_1.$
By \eqref{eq: bound by TV}, we obtain
\begin{equation}
\PP_{\pi}\left(\mathcal{A}_0\right)\geq 1-\alpha_0-\alpha_1-\int \left|f_{\theta_0}\left(z\right)- f_{\pi}\left(z\right)\right| dz.
\label{eq: prob under the alternative set}
\end{equation}
For $z\in \mathcal{A}_0$ and $\theta \in \FF$, by triangle inequality, 
\begin{equation}
\|\widehat{\beta}(z)-\beta(\theta)\|_{q}\geq \left|\|\beta(\theta)-\betanull\|_{q}-\|\widehat{\beta}(z)-\betanull\|_{q}\right| \geq \left(1-c_0\right)\distalter.
\label{eq: dominating term under the alternative}
\end{equation}
For $z\in \mathcal{A}_0$ and $\theta\in \FF$, then
$\left|\widehat{L}_{q}\left(z\right)-\|\widehat{\beta}(z)-\beta\left(\theta\right)\|_{q}^2\right|
\geq  \left|\|\widehat{\beta}(z)-\beta\left(\theta\right)\|_{q}^2-\|\widehat{\beta}(z)-\betanull\|_{q}^2\right|-\left|\widehat{L}_{q}\left(z\right)-\|\widehat{\beta}(z)-\betanull\|_{q}^2\right|
\geq  (1-2c_0-\frac{1}{4}) \distalter^2,$
where the first inequality follows from triangle inequality and the last inequality follows from \eqref{eq: high prob event at null hypothesis} and \eqref{eq: dominating term under the alternative}. Hence, for $z \in \mathcal{A}_0$, we obtain 
\begin{equation}
\inf_{\theta\in \FF}\left|\widehat{L}_{q}\left(z\right)-\|\widehat{\beta}(z)-\beta\left(\theta\right)\|_{q}^2\right|\geq (1-2c_0-\frac{1}{4}) \distalter^2.
\label{eq: cannot estimate alternative well}
\end{equation} 
Note that 
$\sup_{\theta \in \FF} \PP_{\theta} \left(\left|\widehat{L}_{q}\left(Z\right)-\|\widehat{\beta}(Z)-\beta\left(\theta\right)\|_{q}^2\right|\geq (1-2c_0-\frac{1}{4}) \distalter^2\right)
\geq \\ \sup_{\theta \in \FF} \PP_{\theta} \left(\inf_{\theta\in \FF} \left|\widehat{L}_{q}\left(Z\right)-\|\widehat{\beta}(Z)-\beta\left(\theta\right)\|_{q}^2\right|\geq (1-2c_0-\frac{1}{4}) \distalter^2\right).$
Since the max risk is lower bounded by the Bayesian risk, we can further lower bound the last term by $\PP_{\pi} \left(\inf_{\theta\in \FF} \left|\widehat{L}_{q}\left(Z\right)-\|\widehat{\beta}(Z)-\beta\left(\theta\right)\|_{q}^2\right|\geq (1-2c_0-\frac{1}{4}) \distalter^2\right).$ Combined with \eqref{eq: cannot estimate alternative well}, we establish 
\begin{equation}
\sup_{\theta \in \FF} \PP_{\theta} \left(\left|\widehat{L}_{q}\left(Z\right)-\|\widehat{\beta}(Z)-\beta\left(\theta\right)\|_{q}^2\right|\geq (1-2c_0-\frac{1}{4}) \distalter^2\right)\geq \PP_{\pi}(\mathcal{A}_0).
\label{eq: bound over alternative}
\end{equation}
Combining \eqref{eq: bad estimate of risk at theta0 in probability}, \eqref{eq: prob under the alternative set} and \eqref{eq: bound over alternative}, we establish \eqref{eq: key lower bound in high prob unknown}.\\ 
\underline{\textbf{Proof of \eqref{eq: lower bound for CI of loss}}}\\
For $\CIZR\in \II_{\alpha}\left(\Theta, \widehat{\beta},\ell_{q}\right)$, we have
\begin{equation}
\inf_{\theta \in \Theta}\PP_{\theta}\left(\|\widehat{\beta}(Z)-\beta\left(\theta\right)\|_{q}^2 \in \CIZR\right) \geq 1-\alpha.
\label{eq: coverage of CI general}
\end{equation}
Define the event
$\mathcal{A}=\left\{z: \|\widehat{\beta}(z)-\betanull\|_{q}<c_0 \distalter, \; \;  \|\widehat{\beta}(z)-\betanull\|_{q}^2 \in {\rm CI}_{\alpha}\left(\widehat{\beta}, L,z\right)\right\}.$
By \eqref{eq: adaptive risk assumption high prob} and \eqref{eq: coverage of CI general}, we have
$\PP_{\theta_0}\left(\mathcal{A}\right)\geq 1-\alpha-\alpha_0.$
 \eqref{eq: joint prob and marginal prob} and \eqref{eq: bound by TV} imply
\begin{equation}
\pjoint \left(\mathcal{A}\right)=\PP_{\pi}\left(\mathcal{A}\right)\geq 1-\alpha-\alpha_0-\TV\left(f_{\pi},f_{\theta_0}\right).
\label{eq: prob under the alternative general}
\end{equation}
Define the event
$\mathcal{B}_{\theta}=\left\{z: \|\widehat{\beta}(z)-\beta\left(\theta\right)\|_{q}^2 \in {\rm CI}_{\alpha}\left(\widehat{\beta}, \ell_{q},z\right) \right\}$ and $\mathcal{M}=\cup_{\theta\in \FF} \mathcal{B}_{\theta}.$
By \eqref{eq: coverage of CI general}, we have
$$\pjoint \left(\mathcal{M}\right)=\int \left(\int 1_{z\in \mathcal{M} } f_{\theta}(z) dz\right) \pi\left(\theta\right) d\theta \geq \int \left(\int 1_{z \in \mathcal{B}_{\theta} } f_{\theta}(z) dz\right) \pi\left(\theta\right) d\theta \geq 1-\alpha.$$
Combined with  \eqref{eq: prob under the alternative general}, we have 
$\pjoint \left(\mathcal{A}\cap \mathcal{M}\right)\geq 1-2\alpha-\alpha_0-\TV\left(f_{\pi},f_{\theta_0}\right).$
For $z \in \mathcal{M}$, there exists $\bar{\theta}\in \FF$ such that $\|\widehat{\beta}(z)-\beta(\bar{\theta})\|_{q}^2 \in {\rm CI}_{\alpha}\left(\widehat{\beta}, \ell_{q},z\right);$ For $z \in \mathcal{A}$,
we have $\|\widehat{\beta}(z)-\betanull\|_{q}^2 \in {\rm CI}_{\alpha}\left(\widehat{\beta}, \ell_{q},z\right)$ and $\|\widehat{\beta}(z)-\betanull\|_{q}<c_0 \distalter$. Hence, for $z\in \mathcal{A}\cap \mathcal{M}$, we have $\|\widehat{\beta}(z)-\beta(\bar{\theta})\|_{q}^2,\|\widehat{\beta}(z)-\betanull\|_{q}^2 \in {\rm CI}_{\alpha}\left(\widehat{\beta}, \ell_{q},z\right)$  and 
$\|\widehat{\beta}(z)-\beta(\bar{\theta})\|_{q}\geq \|\beta(\bar{\theta})-\betanull\|_{q}-\|\widehat{\beta}(z)-\betanull\|_{q}\geq \left(1-c_0\right)\distalter$
and hence
\begin{equation}
\RR\left({\rm CI}_{\alpha}\left(\widehat{\beta}, \ell_{q},z\right)\right)\geq \left(1-2c_0\right) \distalter^2.
\label{eq: length of CI under alter general}
\end{equation}
Define the event
$\mathcal{C}=\left\{z: \RR\left({\rm CI}_{\alpha}\left(\widehat{\beta}, \ell_{q},z\right)\right)\geq \left(1-2c_0\right)  \distalter^2 \right\}.$
By \eqref{eq: length of CI under alter general}, we have
$\PP_{\pi}\left(\mathcal{C}\right)=\pjoint \left(\mathcal{C}\right)\geq \pjoint \left(\mathcal{A}\cap \mathcal{M}\right) \geq 1-2\alpha-\alpha_0-\TV\left(f_{\pi},f_{\theta_0}\right).$
By \eqref{eq: bound by TV}, we establish 
$\PP_{\theta_0}\left(\mathcal{C}\right) \geq 1-2\alpha-\alpha_0-2\TV\left(f_{\pi},f_{\theta_0}\right)$ and hence \eqref{eq: lower bound for CI of loss}.
\subsection{Proof of Theorems \ref{thm: adaptation of CI for the l2 loss both known} and \ref{thm: adaptivity for q<2 both known}} 
\label{sec: known variance lower bound proof}

We first specify some constants used in the proof. Let $C^{*}$ be given in \eqref{eq: gloss adaptivity assumption}. Define
$\epsilon_1=\frac{1-2\alpha-2\alpha_0}{12}$ and
{\footnotesize
\begin{equation}
c_0=\min\left\{\frac{1}{2},32\log\left(1+\epsilon_1^2\right), \frac{2}{3}\sqrt{\log(1+\epsilon_1^2)}, \frac{1-2\gamma}{16C^{*}},\left(\frac{1-2\gamma}{16C^{*}}\right)^2\right\}, \; c_0^{*}=\min\left\{c_0,\frac{\sqrt{M_1}-1}{C^{*}M_1+\sqrt{M_1}-1}\right\}.
\label{eq: key constant}
\end{equation}
}
Theorems \ref{thm: adaptation of CI for the l2 loss both known} and \ref{thm: adaptivity for q<2 both known} follow from Theorem \ref{thm: general adaptivity for known variance} below. 
\begin{Theorem}
\label{thm: general adaptivity for known variance}
Suppose $0<\alpha<\frac{1}{4}$, $1\leq q \leq 2$ and the sparsity levels $k_1,k_2$ and $k_0$ satisfy Assumption $({\rm B}2)$.  Suppose that $\widehat{\beta}$ satisfies Assumption $({\rm A}2)$ with $\|\betanull\|_0\leq k_0$.
\begin{enumerate}
\item If $k_2 \lesssim \frac{\sqrt{n}}{\log p}$, then there is some constant $c>0$ such that
\begin{equation}
\RR_{\alpha}^{*}\left(\Theta_0\left(k_1\right),\Theta_0\left(k_2\right), \widehat{\beta},\ell_{q}\right) \geq c k_2^{\frac{2}{q}} {\frac{\log p}{n}}\sigma_0^2.
\label{eq: sparse general adaptivity lower bound known variance}
\end{equation}
\item
If $\frac{\sqrt{n}}{\log p}\lesssim k_2 \lesssim \frac{n}{\log p}$, then there is some constant $c>0$ such that 
\begin{equation}
\begin{aligned}
&\RR_{\alpha}^{*}\left(\Theta_0\left(k_1\right),\Theta_0\left(k_2\right), \widehat{\beta},\ell_{q}\right) \\
&\geq c \max\left\{\left((1-c_2)^2 k_2^{\frac{2}{q}-1}k_{1}{\frac{\log p}{n}}-(1+c_1)^2k^{\frac{2}{q}}_{1}{\frac{\log p}{n}}\right)_{+},\frac{k_2^{\frac{2}{q}-1}}{\sqrt{n}}\right\}\sigma_0^2,
\label{eq: dense general adaptivity lower bound known variance}
\end{aligned}
\end{equation}
where $c_1=\frac{C^{*} M_1 k_0^{\frac{1}{q}}}{(k_1-k_0)^{\frac{1}{q}}}$ and $c_2=\frac{C^{*} k_0^{\frac{1}{q}}}{M_1(k_2-k_0)^{\frac{1}{q}-\frac{1}{2}}(k_1-k_0)^{\frac{1}{2}}}$. 
\end{enumerate}
In particular, the minimax lower bound \eqref{eq: sparse general adaptivity lower bound known variance} and the term $\frac{k_2^{\frac{2}{q}-1}}{\sqrt{n}}\sigma_0^2$ in \eqref{eq: dense general adaptivity lower bound known variance} can be established under the weaker assumption $({\rm A}1)$ with $\|\betanull\|_0\leq k_0$.
\end{Theorem}

By Theorem \ref{thm: general adaptivity for known variance}, we establish \eqref{eq: l2 lower bound adaptation} in Theorem \ref{thm: adaptation of CI for the l2 loss both known} and \eqref{eq: lq lower bound adaptation} in Theorem \ref{thm: adaptivity for q<2 both known}. 
In the regime $k_2\lesssim \frac{\sqrt{n}}{\log p}$, the lower bound \eqref{eq: l2 lower bound adaptation} for $q=2$ and \eqref{eq: lq lower bound adaptation} for $1\leq q<2$ follow from \eqref{eq: sparse general adaptivity lower bound known variance}. For the case $q=2$, in the regime $\frac{\sqrt{n}}{\log p}\lesssim k_2 \lesssim \frac{{n}}{\log p}$, the first term of the right hand side of \eqref{eq: dense general adaptivity lower bound known variance} is $0$ while the second term is $\frac{1}{\sqrt{n}}$, which leads to \eqref{eq: l2 lower bound adaptation}. 
For $1\leq q <2$, let $k_1^{*}=\min\{k_1,\zeta_0 k_2\}$ for some constant $0<\zeta_0<1$, an application of \eqref{eq: dense general adaptivity lower bound known variance} leads to  
$\RR_{\alpha}^{*}\left(\Theta_0\left(k_1^{*}\right),\Theta_0\left(k_2\right), \widehat{\beta},\ell_{q}\right) \geq c \max\left\{k_2^{\frac{2}{q}-1}k^{*}_{1}{\frac{\log p}{n}}, \frac{k_2^{\frac{2}{q}-1}}{\sqrt{n}}\right\} \sigma_0^2.$
By this result, if $k_1\leq \zeta_0 k_2$, the lower bounds  \eqref{eq: lq lower bound adaptation} 
in the regions $\regimeb$ and $\regimec$ follow; if $\zeta_0 k_2 <k_1 \leq k_2$, by the fact that $\RR_{\alpha}^{*}\left(\Theta_0\left(k_1\right),\Theta_0\left(k_2\right), \widehat{\beta},\ell_{q}\right)\geq \RR_{\alpha}^{*}\left(\Theta_0\left(k_1^{*}\right),\Theta_0\left(k_2\right), \widehat{\beta},\ell_{q}\right)$ and $k_1^{*}=\zeta_0 k_2 \geq \zeta_0 k_1$, 
the lower bounds \eqref{eq: lq lower bound adaptation} 
over the regions $\regimeb$ and $\regimec$ follow. The following lemma shows that \eqref{eq: l2 lower bound adaptation} holds 
for $\betal$ defined in \eqref{eq: Lasso} with $A>\sqrt{2}$ by verifying Assumption $({\rm A}1)$ and \eqref{eq: lq lower bound adaptation} holds for $\betal$ defined in \eqref{eq: Lasso} with $A>4\sqrt{2}$ by verifying Assumption $({\rm A}2)$. Its proof can be found in the supplement \citep{cai2016asupplement}.
\begin{Lemma}
\label{lem: lemma at null lasso}
If $A>4\sqrt{2}$, then we have 
$$\inf_{\left\{\theta=\left(\betanull ,{\rm I}, \sigma\right):\sigma\leq 2\sigma_0\right\}}\PP_{\theta}\left(\|\betal-\betanull\|_{q}^2\leq C \|\betanull\|_0^{\frac{2}{q}} \frac{\log p}{n}\sigma^2\right)\geq 1-c\exp\left(-c'n\right)-p^{-c}.$$
In particular, the above result holds for $q=2$ under the assumption $A>\sqrt{2}$.
\end{Lemma}

\putbib
\end{bibunit}


\begin{thebibliography}{10}

\bibitem{arias2011global}
Ery Arias-Castro, Emmanuel~J Cand{\`e}s, and Yaniv Plan.
\newblock Global testing under sparse alternatives: Anova, multiple comparisons
  and the higher criticism.
\newblock {\em The Annals of Statistics}, 39(5):2533--2556, 2011.

\bibitem{bayati2012lasso}
Mohsen Bayati and Andrea Montanari.
\newblock The lasso risk for gaussian matrices.
\newblock {\em Information Theory, IEEE Transactions on}, 58(4):1997--2017,
  2012.

\bibitem{belloni2011square}
Alexandre Belloni, Victor Chernozhukov, and Lie Wang.
\newblock Square-root lasso: pivotal recovery of sparse signals via conic
  programming.
\newblock {\em Biometrika}, 98(4):791--806, 2011.

\bibitem{bickel2009simultaneous}
Peter~J Bickel, Ya'acov Ritov, and Alexandre~B Tsybakov.
\newblock Simultaneous analysis of lasso and dantzig selector.
\newblock {\em The Annals of Statistics}, 37(4):1705--1732, 2009.

\bibitem{buhlmann2011statistics}
Peter B{\"u}hlmann and Sara Van De~Geer.
\newblock {\em Statistics for high-dimensional data: methods, theory and
  applications}.
\newblock Springer Science \& Business Media, 2011.

\bibitem{cai2016asupplement}
T~Tony Cai and Zijian Guo.
\newblock Supplement to ``accuracy assessment for high-dimensional linear
  regression".
\newblock 2016.

\bibitem{cai2015regci}
T~Tony Cai and Zijian Guo.
\newblock Confidence intervals for high-dimensional linear regression: Minimax
  rates and adaptivity.
\newblock {\em The Annals of Statistics}, to appear.

\bibitem{cai2004adaptation}
T~Tony Cai and Mark~G Low.
\newblock An adaptation theory for nonparametric confidence intervals.
\newblock {\em The Annals of statistics}, 32(5):1805--1840, 2005.

\bibitem{Cai06}
T~Tony Cai and Mark~G Low.
\newblock Adaptive confidence balls.
\newblock {\em The Annals of Statistics}, 34(1):202--228, 2006.

\bibitem{cai2014adaptive}
T~Tony Cai, Mark~G Low, and Zongming Ma.
\newblock Adaptive confidence bands for nonparametric regression functions.
\newblock {\em Journal of the American Statistical Association},
  109:1054--1070, 2014.

\bibitem{cai2009data}
T.~Tony Cai and Harrison~H Zhou.
\newblock A data-driven block thresholding approach to wavelet estimation.
\newblock {\em The Annals of Statistics}, 37(2):569--595, 2009.

\bibitem{candes2007dantzig}
Emmanuel Cand{\`e}s and Terence Tao.
\newblock The dantzig selector: statistical estimation when p is much larger
  than n.
\newblock {\em The Annals of Statistics}, 35(6):2313--2351, 2007.

\bibitem{chernozhukov2015post}
Victor Chernozhukov, Christian Hansen, and Martin Spindler.
\newblock Post-selection and post-regularization inference in linear models
  with many controls and instruments.
\newblock 2015.

\bibitem{chernozhukov2015valid}
Victor Chernozhukov, Christian Hansen, and Martin Spindler.
\newblock Valid post-selection and post-regularization inference: An
  elementary, general approach.
\newblock {\em arXiv preprint arXiv:1501.03430}, 2015.

\bibitem{donoho1995adapting}
David~L Donoho and Iain~M Johnstone.
\newblock Adapting to unknown smoothness via wavelet shrinkage.
\newblock {\em Journal of the American Statistical Association},
  90(432):1200--1224, 1995.

\bibitem{donoho2011noise}
David~L Donoho, Arian Maleki, and Andrea Montanari.
\newblock The noise-sensitivity phase transition in compressed sensing.
\newblock {\em Information Theory, IEEE Transactions on}, 57(10):6920--6941,
  2011.

\bibitem{guo2016optimal}
Zijian Guo, Wanjie Wang, T~Tony Cai, and Hongzhe Li.
\newblock Optimal estimation of co-heritability in high-dimensional linear
  models.
\newblock {\em arXiv preprint arXiv:1605.07244}, 2016.

\bibitem{hoffmann2011adaptive}
Marc Hoffmann and Richard Nickl.
\newblock On adaptive inference and confidence bands.
\newblock {\em The Annals of Statistics}, 39(5):2383--2409, 2011.

\bibitem{ingster2010detection}
Yuri~I Ingster, Alexandre~B Tsybakov, and Nicolas Verzelen.
\newblock Detection boundary in sparse regression.
\newblock {\em Electronic Journal of Statistics}, 4:1476--1526, 2010.

\bibitem{janson2015eigenprism}
Lucas Janson, Rina~Foygel Barber, and Emmanuel Cand{\`e}s.
\newblock Eigenprism: Inference for high-dimensional signal-to-noise ratios.
\newblock {\em arXiv preprint arXiv:1505.02097}, 2015.

\bibitem{li1985stein}
Ker-Chau Li.
\newblock From stein's unbiased risk estimates to the method of generalized
  cross validation.
\newblock {\em The Annals of Statistics}, 13(4):1352--1377, 1985.

\bibitem{nickl2013confidence}
Richard Nickl and Sara van~de Geer.
\newblock Confidence sets in sparse regression.
\newblock {\em The Annals of Statistics}, 41(6):2852--2876, 2013.

\bibitem{raskutti2011minimax}
Garvesh Raskutti, Martin~J Wainwright, and Bin Yu.
\newblock Minimax rates of estimation for high-dimensional linear regression
  over-balls.
\newblock {\em Information Theory, IEEE Transactions on}, 57(10):6976--6994,
  2011.

\bibitem{robins2006adaptive}
James Robins and Aad Van Der~Vaart.
\newblock Adaptive nonparametric confidence sets.
\newblock {\em The Annals of Statistics}, 34(1):229--253, 2006.

\bibitem{stein1981estimation}
Charles~M Stein.
\newblock Estimation of the mean of a multivariate normal distribution.
\newblock {\em The Annals of Statistics}, 9(6):1135--1151, 1981.

\bibitem{sun2012scaled}
Tingni Sun and Cun-Hui Zhang.
\newblock Scaled sparse linear regression.
\newblock {\em Biometrika}, 101(2):269--284, 2012.

\bibitem{thrampoulidis2015asymptotically}
Christos Thrampoulidis, Ashkan Panahi, and Babak Hassibi.
\newblock Asymptotically exact error analysis for the generalized $
  \ell_2^2$-lasso.
\newblock {\em arXiv preprint arXiv:1502.06287}, 2015.

\bibitem{tibshirani1996regression}
Robert Tibshirani.
\newblock Regression shrinkage and selection via the lasso.
\newblock {\em Journal of the Royal Statistical Society. Series B
  (Methodological)}, 58(1):267--288, 1996.

\bibitem{van2014asymptotically}
Sara van~de Geer, Peter B{\"u}hlmann, YaÕacov Ritov, and Ruben Dezeure.
\newblock On asymptotically optimal confidence regions and tests for
  high-dimensional models.
\newblock {\em The Annals of Statistics}, 42(3):1166--1202, 2014.

\bibitem{verzelen2012minimax}
Nicolas Verzelen.
\newblock Minimax risks for sparse regressions: Ultra-high dimensional
  phenomenons.
\newblock {\em Electronic Journal of Statistics}, 6:38--90, 2012.

\bibitem{ye2010rate}
Fei Ye and Cun-Hui Zhang.
\newblock Rate minimaxity of the lasso and dantzig selector for the $\ell_{q}$
  loss in $\ell_{r}$ balls.
\newblock {\em The Journal of Machine Learning Research}, 11:3519--3540, 2010.

\bibitem{yi2013sure}
Feng Yi and Hui Zou.
\newblock {SURE}-tuned tapering estimation of large covariance matrices.
\newblock {\em Computational Statistics \& Data Analysis}, 58:339--351, 2013.

\end{thebibliography}


\end{document}